\documentclass{article}
	
	\usepackage{subfigure}

\usepackage{graphicx}
\usepackage{url}
\usepackage{pslatex}
\usepackage{bm}
\usepackage{amsmath,amssymb,amsthm}  
\usepackage{paralist}
\usepackage{amsmath}               
{
	\theoremstyle{plain}
	
}
\newtheorem{theorem}{Theorem}

\newtheorem{lemma}{Lemma}

\newtheorem{remark}{Remark}
\usepackage{algorithm}
\usepackage{algorithmic}        
\usepackage{multirow}

\renewcommand{\algorithmicrequire}{\textbf{Input:}} 
\renewcommand{\algorithmicensure}{\textbf{Iteration:}}
\renewcommand{\algorithmiclastcon}{\textbf{Output:}} 

\newcommand{\e}{\begin{equation}}
\newcommand{\ee}{\end{equation}}
\newcommand{\en}{\begin{equation*}}
\newcommand{\een}{\end{equation*}}
\newcommand{\eqn}{\begin{eqnarray}}
\newcommand{\eeqn}{\end{eqnarray}}
\newcommand{\bmat}{\begin{bmatrix}}
\newcommand{\emat}{\end{bmatrix}}
\newcommand{\BIT}{\begin{itemize}}
\newcommand{\EIT}{\end{itemize}}

\newcommand{\argmin}{\mathop{\rm argmin}}

\newcommand{\ve}{\bm e}
\newcommand{\vf}{\bm f}

\newcommand{\vr}{\bm r}

\newcommand{\vu}{\bm u}
\newcommand{\vv}{\bm v}

\newcommand{\vx}{\bm x}
\newcommand{\vy}{\bm y}

\newcommand{\vsigma}{\bm \sigma}

\newcommand{\mA}{\bm A}
\newcommand{\mB}{\bm B}

\newcommand{\mE}{\bm E}

\newcommand{\mI}{\bm I}

\newcommand{\mGamma}{\bm \Gamma}

\newcounter{oursection}

	\usepackage{hyperref}
	\usepackage{cleveref}
	\usepackage{appendix}
	\usepackage{soul}
	\usepackage{mathtools, cuted,multicol}
	\crefname{assumption}{assumption}{assumptions}
	\usepackage{arxiv}
	\usepackage{cases}
	\usepackage[utf8]{inputenc} 
	\usepackage[T1]{fontenc}    
	\usepackage{url}            
	\usepackage{booktabs}       
	\usepackage{amsfonts}       
	\usepackage{nicefrac}       
	\usepackage{microtype}      
	\usepackage{lipsum}
	\usepackage{color}
	\usepackage{tikz}
	\usepackage{pgfplots}
	\pgfplotsset{compat=1.5.1}
	\usepgfplotslibrary{groupplots}
	\usetikzlibrary{intersections,calc,arrows,matrix,spy,math,calc,fit,positioning}

	\usetikzlibrary{arrows.meta,shapes.geometric}
	
	\newcommand*\circled[1]{\tikz[baseline=(char.base)]{
			\node[shape=circle,draw,inner sep=0.05pt] (char) {#1};}}
	\title{On Adapting Nesterov's Scheme to Accelerate Iterative Methods for Linear Problems}
	
	\author{
	 Tao Hong\\
	Technion - Israel Institute of Technology\\
	\texttt{hongtao@cs.technion.ac.il}
	\And
	  Irad Yavneh \\
	Technion - Israel Institute of Technology\\
	\texttt{irad@cs.technion.ac.il}
	}
	
\begin{document}
	
	\maketitle
	
	\begin{abstract}
	Nesterov's well-known scheme for accelerating gradient descent in convex optimization problems is adapted to accelerating stationary iterative solvers for linear systems. Compared with classical Krylov subspace acceleration methods, the proposed scheme requires more iterations, but it is trivial to implement and retains essentially the same computational cost as the unaccelerated method. An explicit formula for a fixed optimal parameter is derived in the case where the stationary iteration matrix has only real eigenvalues, based only on the smallest and largest eigenvalues. The fixed parameter, and corresponding convergence factor, are shown to maintain their optimality when the iteration matrix also has complex eigenvalues that are contained within an explicitly defined disk in the complex plane. A comparison to Chebyshev acceleration based on the same information of the smallest and largest real eigenvalues (dubbed Restricted Information Chebyshev acceleration) demonstrates that Nesterov's scheme is more robust in the sense that it remains optimal over a larger domain when the iteration matrix does have some complex eigenvalues.  Numerical tests validate the efficiency of the proposed scheme.  This work generalizes and extends the results of \cite[Lemmas 3.1 and 3.2 and Theorem 3.3]{de2021asymptotic}.
	
	\end{abstract}
	
	

\section{Introduction}
Many scientific computing applications require solving linear systems of equations of the form \cite{brandt1977multi,saad2003iterative}: 
	\begin{equation}
	\mA\vx=\vf,\label{eq:linearPro}
	\end{equation}
	where $\mA\in\mathbb R^{N\times N}$ is a sparse, large-scale, ill-conditioned matrix. For example, $\mA$ may be a  discretization of an elliptic partial differential equation (PDE) or system. Because direct solvers are relatively expensive, especially for 3D problems, iterative methods are often preferred, e.g., successive over-relaxation or multigrid. These are very often used to advantage as preconditioners for Krylov subspace acceleration methods. {\color{black}The LOPCG method for eigenvalue problems\cite{knyazev2001toward,knyazev2003efficient} is an alternative acceleration method, which uses a linear combination of two consecutive iterates, together with a preconditioned residual, to construct the next iterate such that the residual norm at the current step is minimized.} Motivated by Nesterov's scheme developed in the framework of convex optimization, we consider adapting this approach to accelerating iterative methods for linear problems, using a fixed optimal scalar parameter for which we derive an explicit formula.  
	
	Nesterov's well-known scheme for accelerating gradient descent for convex optimization problems \cite{nesterov1983method,nesterov2018lectures} has attracted much attention in the optimization community over the years. Given an unconstrained optimization problem,
	$$
	\vx^*=\arg\min_{\vx\in \mathbb R^N}  F(\vx),
	$$
	the $k$th iteration of Nesterov's scheme is defined as follows: 
	\begin{equation}
	\begin{array}{rcl}
	\vx_{x+1}&=&\mB(\vy_k),\\
	\vy_{k+1}&=&\vx_{k+1}+c_k(\vx_{k+1}-\vx_k),
	\end{array}\label{eq:Nesterov}
	\end{equation}
	where $\vy_k\in\mathbb R^N$ is an intermediate iterate and $\mB(\cdot)$ represents some iterative method, e.g., gradient descent. In the classical method, an optimal analytical sequence $c_k$ is introduced for convex problems with $\mB(\cdot)$ defined by gradient descent \cite{nesterov2018lectures}. We refer the reader to \cite{su2014differential,wibisono2016variational} and references therein for further discussion of Nesterov's scheme and the choice of $c_k$ {\color{black} with some restarting strategies.} Over the years, Nesterov's scheme has been extended to more general $\mB(\cdot)$ operators, including the proximal operator \cite{beck2009fast}, coordinate descent \cite{nesterov2017efficiency}, {\color{black}alternating least squares \cite{mitchell2020nesterov}}, second-order methods  \cite{nesterov2008accelerating}, and stochastic methods \cite{nitanda2014stochastic}.
	
	In this paper, we focus on linear systems \eqref{eq:linearPro} and adapt \eqref{eq:Nesterov} to acceleration of stationary iterative methods. Since the problem of interest is linear, we simplify \eqref{eq:Nesterov} as follows:
	\begin{equation}
	\begin{array}{rcl}
	\vx_{x+1}&=&\mB\vy_k+\text{Constant},\\
	\vy_{k+1}&=&\vx_{k+1}+c_k(\vx_{k+1}-\vx_k),
	\end{array}\label{eq:NesterovLinear:General}
	\end{equation}
	where we use the fixed matrix $\mB\in \mathbb R^{N\times N}$ (called iteration matrix) plus an appropriate constant vector to represent the operator $\mB(\cdot)$. Below, we show that if all the eigenvalues of $\mB$ are real, then an optimal fixed $c_k$ can be found analytically, depending on the smallest and largest eigenvalues of $\mB$. We note that a similar analysis was performed in \cite{de2021asymptotic} for the restricted case of gradient descent and the assumption that all eigenvalues of $\mB$ are real and of the same sign. {\color{black}Furthermore, \cite{de2021asymptotic} also studied the complex eigenvalues case and proved a lower and upper ACF bound when all complex eigenvalues lie in a prescribed rectangular region.} Since the additional computations in \eqref{eq:NesterovLinear:General} are negligible, we suggest that \eqref{eq:NesterovLinear:General} may in some cases prove to be competitive for accelerating stationary iterative methods. Because $c_k$ will be fixed, we discard the subscript $k$ in the rest of the paper. Interestingly, we find that the results developed for $\mB$ whose eigenvalues are real are also valid for $\mB$ which has complex eigenvalues within a relatively large domain. Moreover, a ``valid'' region in the  complex plane, defined as a region where existence of complex eigenvalues of $\mB$ does not influence the optimal $c$ or the asymptotic convergence factor (ACF), is explicitly identified, dependent on the smallest and largest real eigenvalues. Furthermore, we compare Nesterov's scheme to a ``restricted-information'' (RI) Chebyshev acceleration, where we choose the optimal parameters based on the same information as required for Nesterov's scheme, i.e., the smallest and largest real eigenvalues of $\mB$ \cite{hageman2012applied}. Our comparison indicates that Nesterov's scheme is more robust than RI Chebyshev acceleration with respect to the existence of complex eigenvalues for $\mB$.
	
	The rest of this paper is organized as follows. The analytical derivation of $c$ for $\mB$ with only real eigenvalues is given in \Cref{sec:MainResultsfixedc_k}. The robustness of these results in cases where $\mB$ also has complex eigenvalues is studied in \Cref{sec:ComplexEigenvalues}, including a comparison to RI Chebyshev acceleration. In \Cref{sec:NumExps} we demonstrate the usefulness of the proposed method in accelerating multigrid solution of some second-order elliptic boundary value problems, and conclusions and future work are summarized in \Cref{sec:ConclusionandFuture}.

\section{Optimal Acceleration}\label{sec:MainResultsfixedc_k}
	Representing $\vy_k$ by the previous iterates $\vx_k$ and $\vx_{k-1}$, we obtain
	\begin{equation}
	\vx_{k+1} = \mB\left(\vx_{k}+c\left(\vx_k-\vx_{k-1}\right)\right)+\text{Constant}.\label{eq:updateAll:x}
	\end{equation}
	Denote $\ve_k=\vx_k-\vx^*$,  where $\vx^*$ is the sought solution. Subtracting $\vx^*$ from both sides of \eqref{eq:updateAll:x} and substituting $\mB\vx^*+\text{Constant}=\vx^*$, we arrive at 
	\begin{equation}
	\ve_{k+1} = \mB\left( \left(1+c\right)\ve_k-c\ve_{k-1} \right).\label{eq:updateError}
	\end{equation}
	Denote $\mE_k = \bmat \ve_k \\ \ve_{k-1}\emat$ and rewrite \eqref{eq:updateError} as 
	$
	\mE_k = \bmat (1+c)\mB & -c\mB \\ \mI & \bm 0\emat \mE_{k-1}.
	$
	Then, the asymptotic convergence factor ACF of \eqref{eq:NesterovLinear:General} is given by $\rho(\bm\Gamma)$, the spectral radius of 
	$
	\bm\Gamma \triangleq  \bmat (1+c)\mB & -c\mB \\ \mI & \bm 0\emat.
	$
	Evidently, if there is a $c$ yielding $\rho(\bm\Gamma)<\rho(\mB)$, then \eqref{eq:NesterovLinear:General} provides acceleration. The following derivation produces a $c$ which minimizes $\rho(\bm\Gamma)$, and this optimal $c$ can easily be calculated analytically if all the eigenvalues of $\mB$ are real and its smallest and largest eigenvalues are known. We note that the cost of obtaining the required information of the smallest and largest eigenvalues of $\mB$ may be low for certain specific linear problems and stationary iterative methods, as discussed in \Cref{sec:NumExps}.
	
	Denote by $\lambda$ and $\bmat \vv_1 \\ \vv_2 \emat$ an eigenvalue and corresponding eigenvector of $\mGamma$. We have
	$
	\bmat (1+c)\mB & -c\mB \\ \mI & \bm 0\emat \bmat\vv_1 \\ \vv_2 \emat = \lambda \bmat\vv_1 \\ \vv_2 \emat,
	$
	yielding
	\begin{equation}
	\left(1+c-\frac{c}{\lambda}\right)\mB\vv_1 = \lambda\vv_1 \label{eq:calculateEigGamma}, 
	\end{equation}
	where the trivial case, $\lambda=0$, is omitted. Let $b\in\mathbb R$ denote some eigenvalue of $\mB$ corresponding to $\lambda$.  From \eqref{eq:calculateEigGamma}, we have $\left(1+c-\frac{c}{\lambda}\right)b = \lambda$, yielding $\lambda^2 - (1+c)b\lambda + cb=0$, with solutions $\lambda_{1}(c,b) = \frac{1}{2}\left[ (1+c)b + \sqrt{(1+c)^2b^2-4cb}\right]$, and $\lambda_{2}(c,b) = \frac{1}{2}\left[ (1+c)b - \sqrt{(1+c)^2b^2-4cb}\right]$. We use 
	\begin{equation}
	 b_{cr} = \frac{4c}{(1+c)^2}
	 \label{eq:bcr}
	\end{equation} 
	to denote the ``critical'' value of $b$ for which the square root term in $\lambda_{1,2}$ vanishes for a given $c$.
	\begin{remark}\label{rem:bcr}
	$\lambda_{1,2}$ are complex if and only if $b$ and $c$ are of the same sign and $0 < |b| < |b_{cr}| \, .$
	\end{remark}
	
	Denote 
	\begin{equation}
	\label{eq:r}
	r(c,b) \triangleq \max\left\{\left|\lambda_{1}(c,b)\right|,\left|\lambda_{2}(c,b)\right|\right\}.
	\end{equation}
	
	Without loss of generality, given that {\color{black} $\mB$  represents the iteration matrix of a convergent method}   with real eigenvalues, we assume that its eigenvalues are ordered as $-1<b_1\leq b_2\leq \cdots \leq b_N<1$. The optimal $c$, given by $c^*=\arg\min_{c}\max_b r(c,b)$, is derived next. To simplify the presentation, we assume $b \neq 0$ throughout, since in this case  $r(c,0)=0$, so a vanishing $b$ has no influence on $c^*$.

	\begin{lemma}\label{lemma:regionC}
		$|c^*|< 1$.
	\end{lemma}
	\begin{proof}
		See \Cref{app:sec:proof:lemma:regionC}.
	\end{proof}
	In light of \Cref{lemma:regionC}, we henceforth restrict $|c|$ to be smaller than 1. This and the fact that $b$ is real imply
	\begin{equation}
	r(c,b) = \frac{1}{2} \left|(1+c)b +sgn(b) \sqrt{(1+c)^2b^2-4cb}\right|,
	\label{eq:expression:r(c,b)}
	\end{equation}
	where $sgn(\cdot)$ is the sign function.  
	
	\begin{remark} \label{rem:imaglambda}
	When the square root term is imaginary (see \Cref{rem:bcr}), we obtain  $r(c,b) = |\lambda_1| = |\lambda_2| = \sqrt{cb}$. 
	\end{remark}
	
	For convenience, we henceforth use $r_c(b)$ (respectively, $r_b(c)$) to denote $r(c,b)$ considered as a single-variable function with a fixed $c$ (respectively, fixed $b$).  First we show that the maximal $r_c(b)$ depends only on the extreme eigenvalues of $\mB$.
	
	\bigskip
	
	\begin{lemma}
	   \label{lemma:extreme_rb}
	   For $|c| < 1$, $r_c(b)$ is maximized at either $b_1$ or $b_N$ and has no local maximum.
	\end{lemma}
	\begin{proof}
	See \Cref{app:sec:proof:lemma:extreme_rb}.
	\end{proof}
	
	Next, we fix $b$ and minimize $r_b(c)$. To this end we first identify the range of $c$ for which $\lambda_{1,2}$ are real. 
	\begin{remark}
	\label{rem:ccr}
	Define the critical $c \in (-1,1)$, {\color{black}for which the square-root term in  \eqref{eq:expression:r(c,b)} vanishes}, by
	$$
	c_{cr}(b)=\frac{1-\sqrt{1-b}}{1+\sqrt{1-b}} \, .
	$$
	Then, {\color{black} by the solution of \eqref{eq:calculateEigGamma},} $\lambda_{1,2}$ are complex if and only if
	$\{ c > c_{cr}(b)$ and $b > 0 \}$ 
	or 
	$\{c < c_{cr}(b)$ and $b < 0 \}$ {\color{black}.}
	\end{remark}
	
	We note that $c_{cr}(b)$ is continuous ($\lim_{b\rightarrow 0}c_{cr}(b)=0$) and monotonically increasing on $b\in(-1,1)$, with $c_{cr}(b)\in (-3+2\sqrt{2},1)$. We next show that $c_{cr}$ is the optimal value of $c$ for minimizing $r_b(c)$.
	\begin{lemma}
		\label{lemma:extreme_rc}
		$$
		  c_{cr}(b) =\argmin_c r_b(c).
		$$
	\end{lemma}
	\begin{proof}
		See \Cref{app:sec:proof:lemma:extreme_rc}.
	\end{proof}

	Summarizing, we observe that, by \Cref{lemma:extreme_rb}, $c^*$ depends only on $b_1$ and $b_N$, and therefore, by \Cref{lemma:extreme_rc}, there are only three possible values of $c$ to consider as candidates for $c^*$: $c_{cr}(b_1)$, $c_{cr}(b_N)$, and the value of $c$ which minimizes $r(c,b)$ subject to $r(c,b_1) = r(c,b_N)$. We map out the regions where each of these three options yields the optimal $c$ in the following theorem.  
	
	\begin{theorem}
	\label{them:optimalC}	
	Let $-1< b_1\leq b_N<1$. Then, the optimal coefficient $c^*$ is given by 
	$$
	c^*=c_{cr}(g(b_1,b_N)),
	$$
	where 
	$$
	g(b_1,b_N)=\left\{ 
	\begin{array}{cl}
	b_N,&b_N\geq -3b_1,\\
	-\frac{8b_Nb_1(b_1+b_N)}{(b_1-b_N)^2},&-\frac{1}{3}b_1<b_N<-3b_1,\\
	b_1,& b_N\leq -\frac{1}{3}b_1,
	\end{array}
	\right.
	$$
	yielding the corresponding asymptotic convergence factor
	$$
	r^*=\left\{ 
	\begin{array}{cl}
	1-\sqrt{1-b_N},&b_N\geq -3b_1,\\
	r(c^*,b_1)=r(c^*,b_N),&-\frac{1}{3}b_1<b_N<-3b_1,\\
	\sqrt{1-b_1}-1,& b_N\leq -\frac{1}{3}b_1.
	\end{array}
	\right.
	$$
	\end{theorem}
	\begin{proof}
	See \Cref{app:sec:proof:theorem:them:optimalC}.
	\end{proof}

	
	
\begin{figure}[!htb]
	    \centering
	    \subfigure[]{
	 \tikzstyle{every node}=[font=\large]
	\begin{tikzpicture}
	\begin{axis}[width=3in,grid=both,grid style={line width=.1pt, draw=gray!10},clip=false,grid=major,
	    x=0.15cm, y=0.15cm,
	    axis lines=middle,
	    xmin=-25,xmax=25,
	    ymin=-25,ymax=25,
	    xtick={-20,-10,0,10,20},
	    xticklabels={-1,-0.5,0,0.5,1},
	    ytick={-20,-10,0,10,20},
	    yticklabels={-1,-0.5, 0,0.5,1},
	    inner axis line style={=>},
	    xlabel=$b_1$,
	    ylabel=$b_N$,
	    xlabel style={at={(ticklabel* cs:1)},anchor=north west,font=\large},
	    ylabel style={at={(ticklabel* cs:1)},anchor=south west,font=\large}
	]
	 \addplot[samples=11,domain=-20/3:0,black, ultra thick] {-3*x};
	 \addplot[samples=11,domain=-20:0,black, ultra thick] {-1/3*x};
	 \addplot[samples=11,domain=-20:20,black, ultra thick]{x};
	 \addplot[samples=11,domain=-20:20,black, ultra thick]{20};
	 \addplot[black, ultra thick] coordinates {
	(-20,-20)
	(-20,-15)
	(-20,-10)
	(-20,0)
	(-20,10)
	(-20,15)
	(-20,20)
	};
	\node[anchor=east] at (axis cs: 8,13) {$T_{top}$};
	\node[anchor=east] at (axis cs: -10,13) {$T_{mid}$};
	\node[anchor=east] at (axis cs: -12,-6) {$T_{bot}$};
	
	
	\end{axis}
	\end{tikzpicture}\label{fig:domain_c_opt:a}
	}
 \subfigure[]{\includegraphics[width=7.5cm,height=7cm]{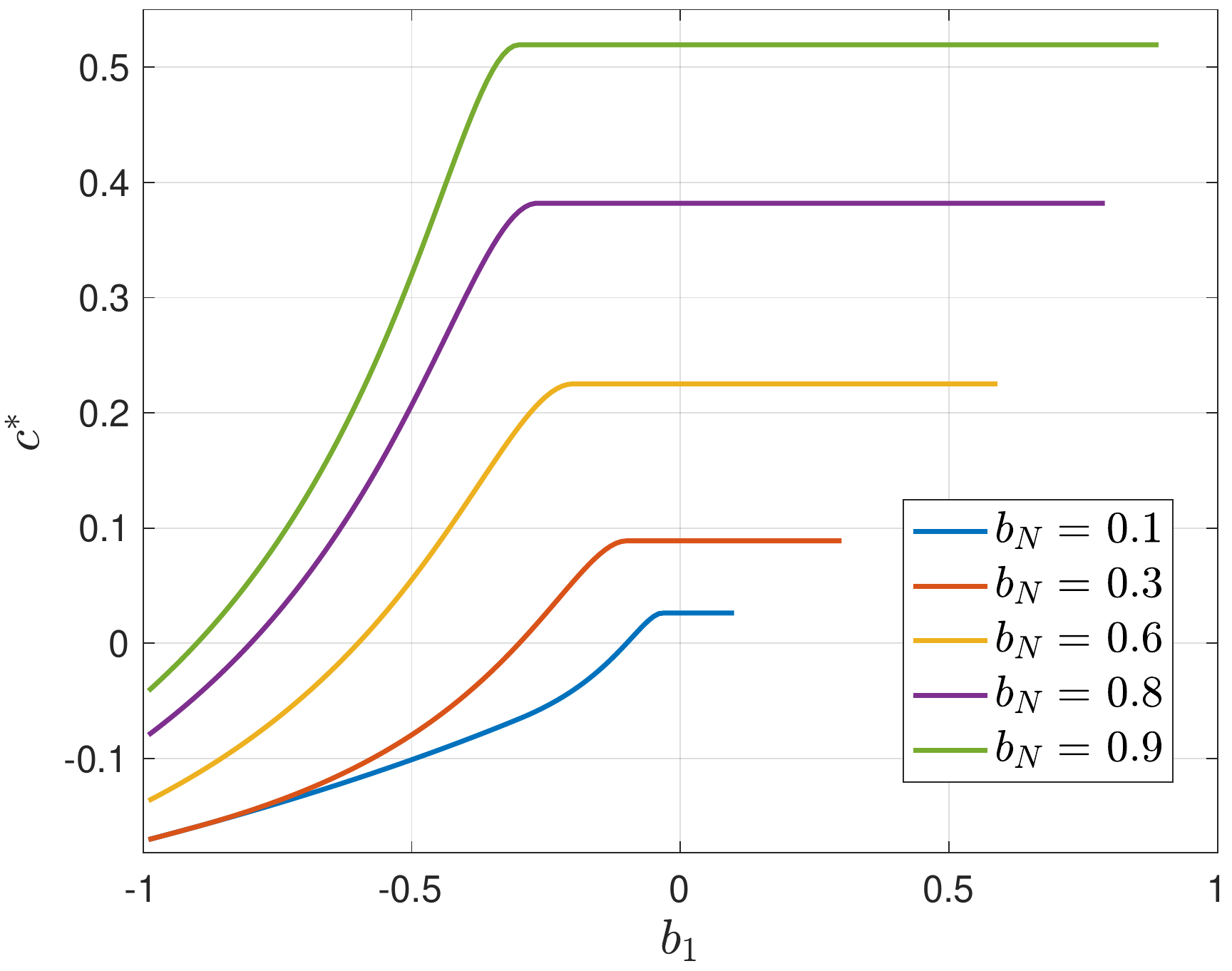}\label{fig:domain_c_opt:b}}
	
	\caption{(a): Three regimes for determining $c^*$ and $r^*$. (b): {\color{black} The value of $c^*$ as a function of $b_1$ and $b_N$.}}
	\label{fig:domain_c_opt}
	\end{figure}
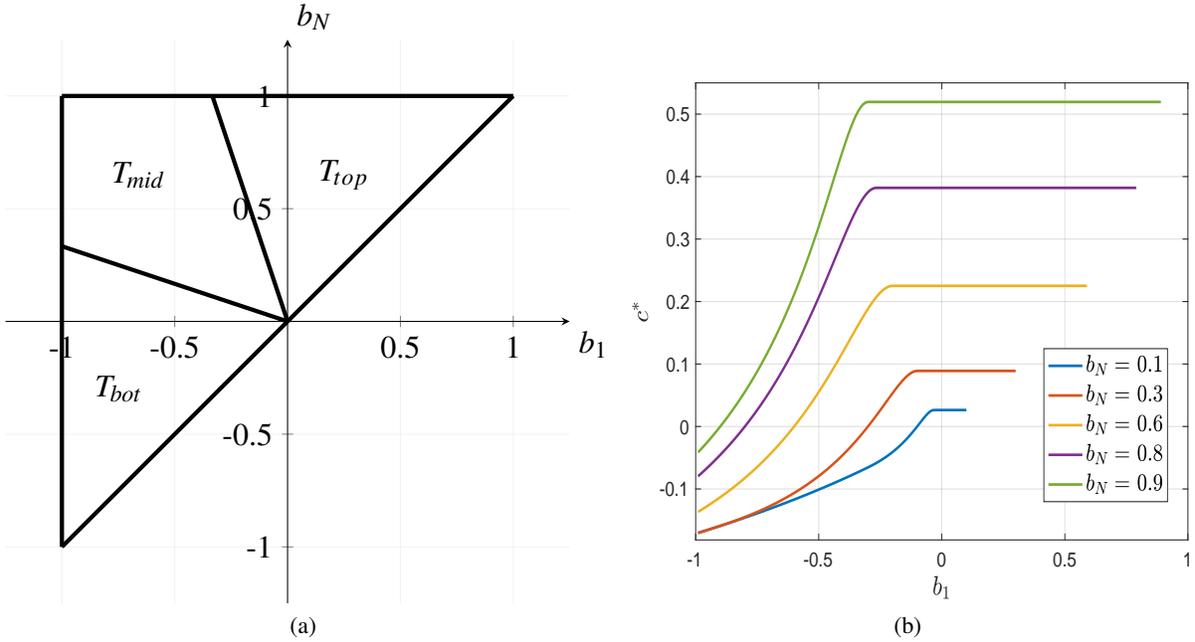

	\begin{remark}
	\label{remark:divergentB}
		{\color{black} Nesterov's scheme can converge even for some $\mB$ whose spectral radii are larger than $1$. In our setting, we can relax our assumption from $-1<b_1\leq b_N<1$ to $-3<b_1\leq b_N<1$, and \Cref{them:Robustness,them:optimalC} remain valid.}
	\end{remark}
Henceforth use $c_{top}$, $c_{mid}$, and $c_{bot}$, as defined in \Cref{app:sec:proof:theorem:them:optimalC}, to denote the optimal coefficient $c^*$ corresponding to regime $T_{top}$ {\color{black}($b_N\geq-3b_1$)}, $T_{mid}$ {\color{black}($b_N\in(-\frac{1}{3}b_1,-3b_1)$)} and $T_{bot}$ {\color{black}($b_N\leq-\frac{1}{3}b_1$)} of \Cref{them:optimalC}, respectively, (see \Cref{fig:domain_c_opt:a}). {\color{black}Also, we numerically show the value of $c^*$ as a function of $b_1$ and $b_N$ in \Cref{fig:domain_c_opt:b}. The flat parts of the curves in \Cref{fig:domain_c_opt:b} imply that $c^*$ only depends on $b_N$ in this regime, and the coincidence of the curves corresponding to $b_N=0.1$ and $b_N=0.3$ for $b_1\rightarrow-1$ indicates that $c^*$ only depends on $b_1$ in that regime.}

This concludes our derivation for the case where the eigenvalues of $\mB$ are all real. In the next section, we extend our results to certain cases where some of the eigenvalues of $\mB$ are complex.
	
\begin{figure}[!htb]
	    \centering
	    \subfigure[]{\includegraphics[scale=0.55]{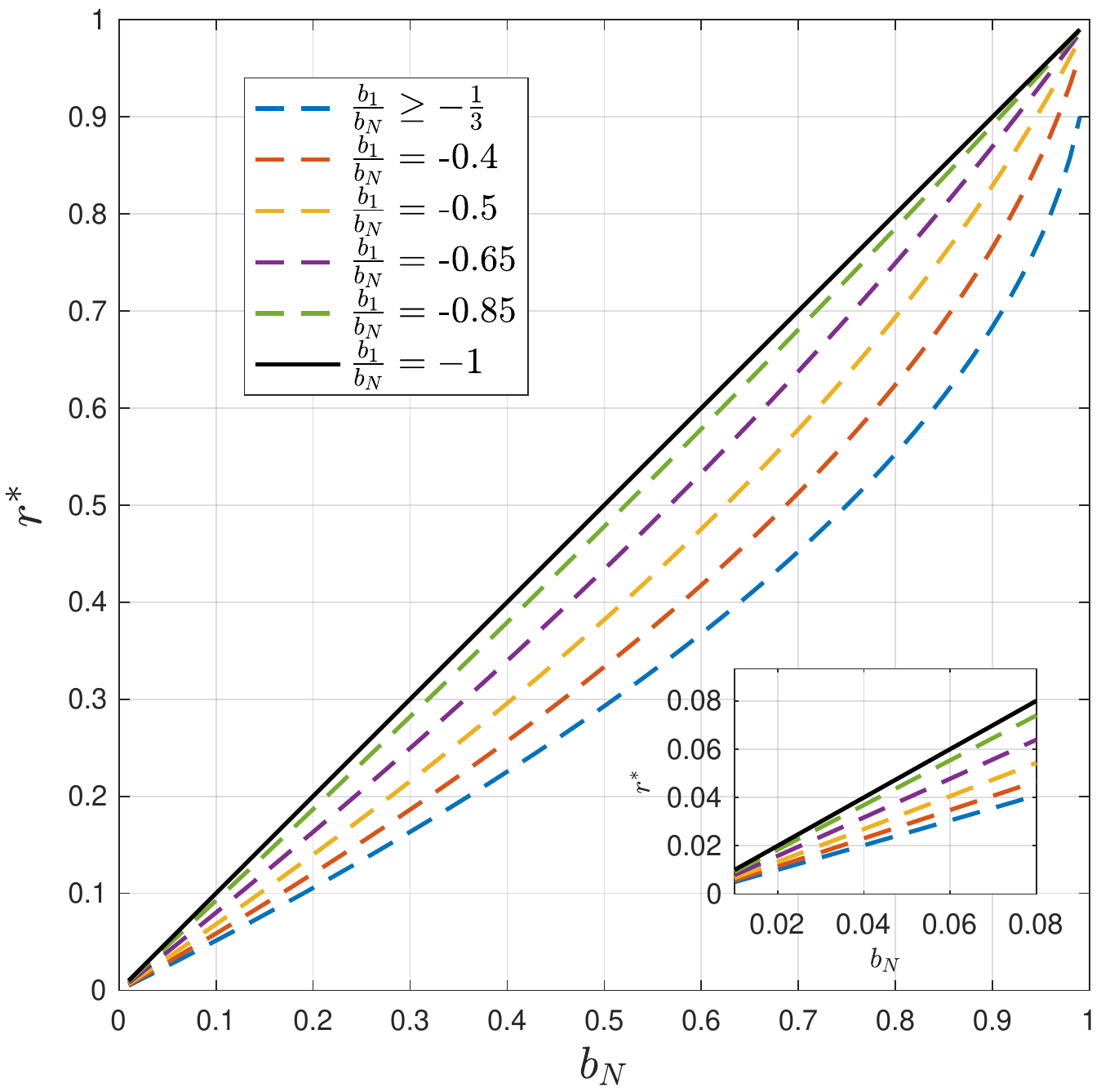}\label{fig:ConvRateNumericallyCheck}}
	    \hspace{0.2cm}
	    \subfigure[]{\includegraphics[scale=0.55]{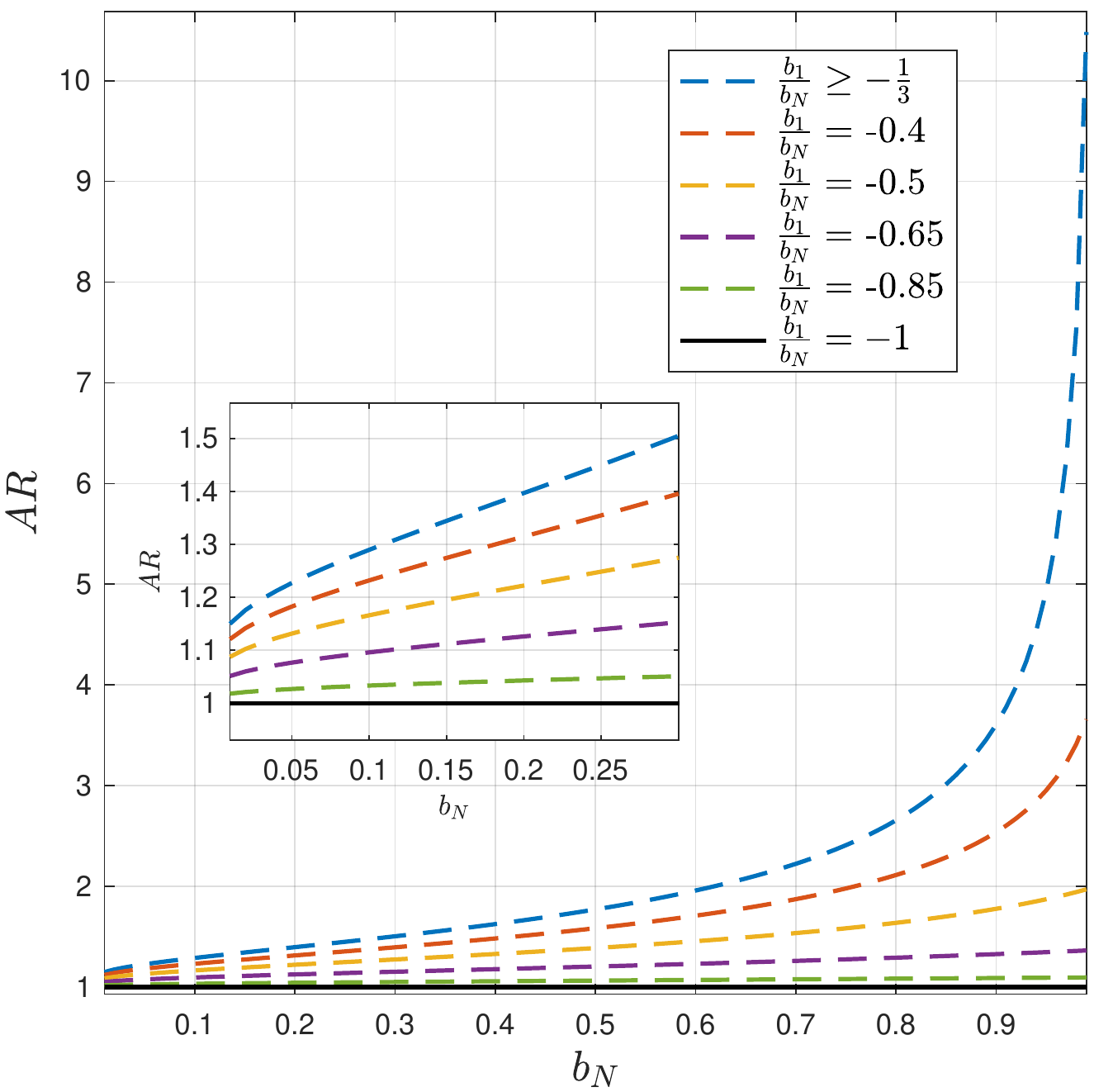}\label{fig:SaveNumericallyCheck}}
	    \caption{(a): The ACF achieved by Nesterov's scheme as a function of $b_N$ and the ratio $\frac{b_1}{b_N}$. (b): The acceleration ratio $AR$ (cf. \eqref{eq:CompSave}) as a function of $b_N$ and the ratio $\frac{b_1}{b_N}$.}
	    \label{fig:NumericallyCheck}
\end{figure}

To conclude this section, we numerically evaluate the savings in computations that are provided by employing Nesterov's scheme \eqref{eq:NesterovLinear:General}. Without loss of generality, we assume $b_N \geq -b_1$. Then, the Acceleration Ratio (AR), defined as the ratio of the number of iterations required without acceleration to the number of iterations required with acceleration, in order to reach the same accuracy, is asymptotically given by 
\begin{equation}
AR = \frac{\log r^*}{\log b_N}\label{eq:CompSave} \, .
\end{equation}
\Cref{fig:NumericallyCheck} shows the relation between $r^*$ and $b_N$ in panel (a), and the acceleration ratio AR in panel (b), for a range of $b_1$ and $b_N$ values. In each curve we fix the ratio $b_1 / b_N$ and vary $b_N\in(0,1)$. For $b_1 / b_N = -1$, there is no acceleration, as discussed above, but for larger ratios we observe acceleration which improves as $b_1 / b_N$ increases, and rapidly improves as $b_N \rightarrow 1$ (that is, when the unaccelerated iterations converge slowly). {\color{black}Note that all values of $r^*$ and $AR$ are identical for $b_1/b_N\geq -1/3$.} 

\section{Complex Eigenvalues}\label{sec:ComplexEigenvalues}
	
	
	\Cref{them:optimalC} is formulated under the assumption that all the eigenvalues of $\mB$ are real. We next state sufficient conditions under which \Cref{them:optimalC} continues to hold even though some of the eigenvalues of $\mB$ are complex.   
	
	\begin{theorem} \label{them:Robustness}
	Assume that, in addition to the real eigenvalues $-1 < b_1 \leq \ldots \leq b_N < 1$ of $\mB$ as assumed in \Cref{them:optimalC}, $\mB$ also has complex eigenvalues. Denote the complex eigenvalues generically by $b^c = \bar b^c e^{j\theta}$, where $j$ is the imaginary unit, $\bar b^c$ is the modulus, and $\theta \in(-\pi,\pi]$ is the argument. Then, $c^*$ and $r^*$ of \Cref{them:optimalC} remain valid if for every one of the complex eigenvalues of $\mB$ the modulus satisfies
	\begin{equation*}
	\bar b^c \leq \left\{\begin{array}{ll}
	\frac{1}{3}b_N     &  c^* = c_{top} \\
	\min(|b_1|,|b_N|)     & c^* = c_{mid} \, .\\
	-\frac{1}{3}b_1  & c^* = c_{bot} 
	\end{array}\right.
	\end{equation*}
	\end{theorem}
	
	%
	
	\begin{proof}
		See \Cref{app:sec:proof:theorem:OptcComplexEigenvaluesB}.
	\end{proof}

\begin{figure*}
		\centering
		\subfigure[$b_1=-0.3$ and $b_N=0.9$.]{ \includegraphics[scale=0.55]{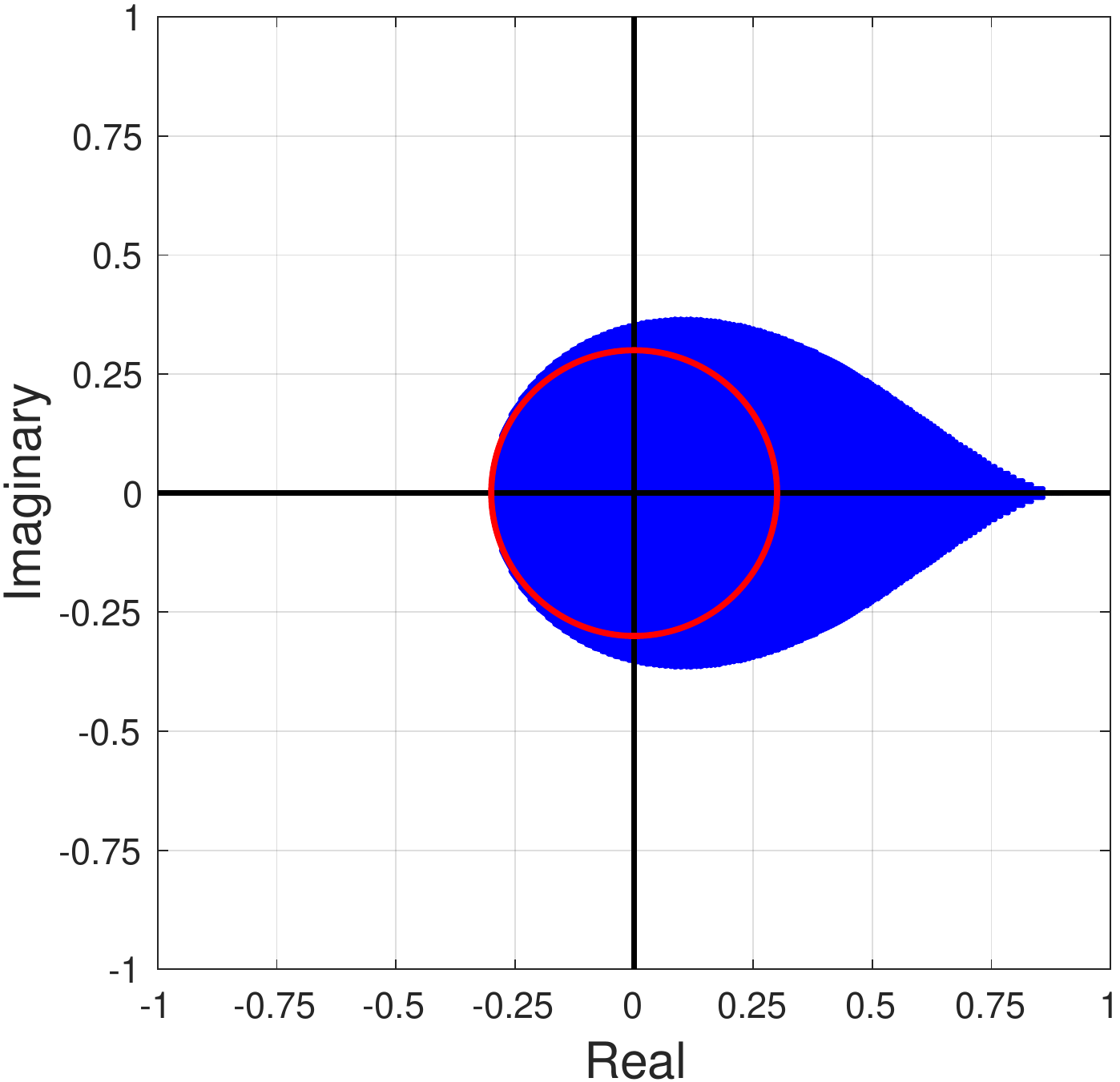}\label{fig:ComplexEigRegions:c1First}}
		\hspace{0.2cm}
		\subfigure[$b_1=-0.5$ and $b_N=0.9$.]{ \includegraphics[scale=0.55]{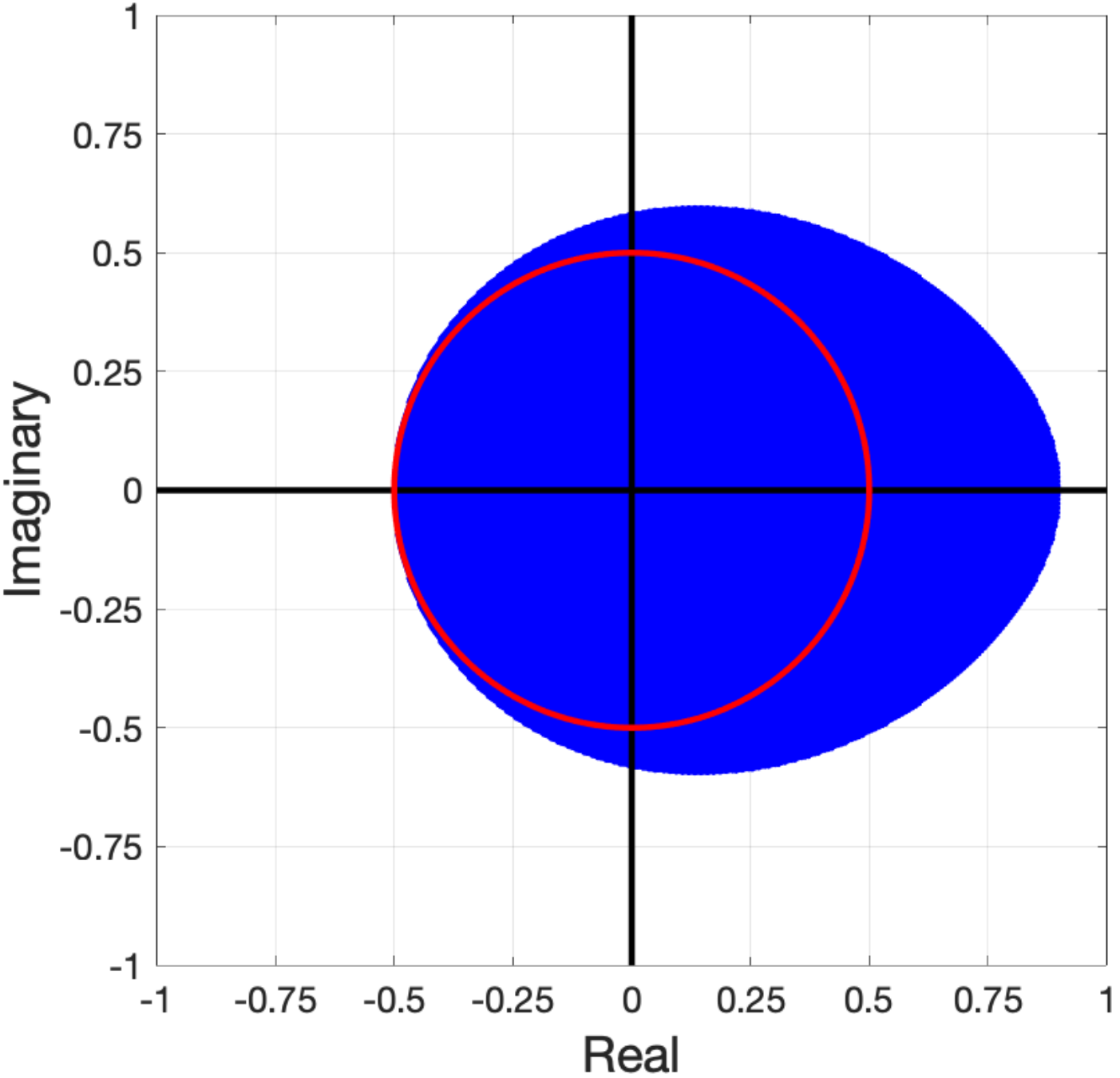}\label{fig:ComplexEigRegions:c2First}}
	
		\subfigure[$b_1=-0.9$ and $b_N=0.3$.]{ \includegraphics[scale=0.55]{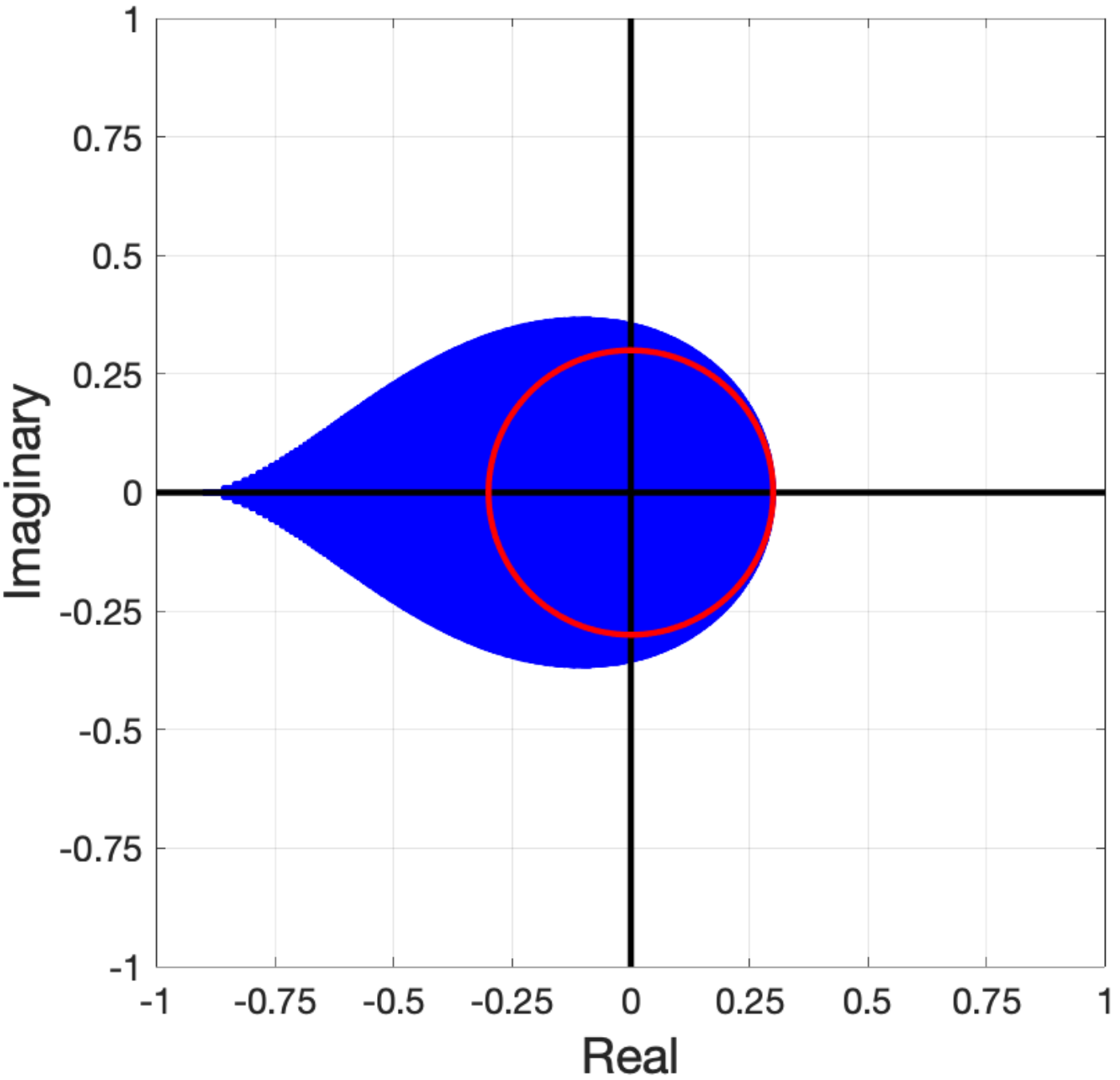}\label{fig:ComplexEigRegions:c3second}}
		\hspace{0.2cm}
		\subfigure[$b_1=-0.9$ and $b_N=0.5$.]{ \includegraphics[scale=0.55]{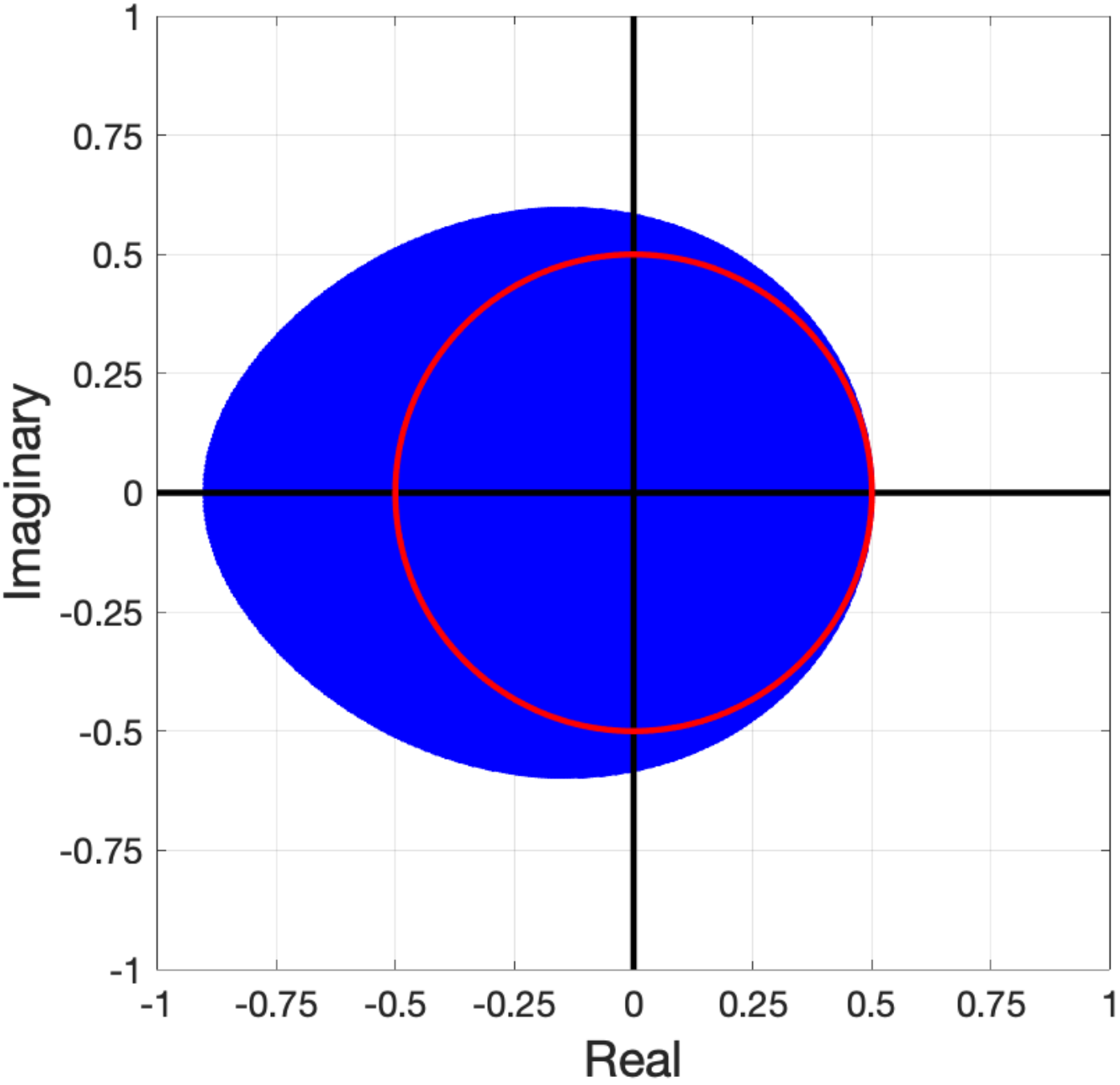}\label{fig:ComplexEigRegions:c2second}}
		\caption{The complex domains defined in \Cref{them:Robustness}. The red circle has radius $|b_1|$ in panels (a) and (b) and $|b_N|$ in (c) and (d).}
		\label{fig:ComplexEigRegions}
	\end{figure*}

	The sufficient robustness conditions of \Cref{them:Robustness} are tight for $|\theta| \rightarrow \pi$ when $b_N > |b_1|$, and can be relaxed increasingly as $|\theta|$ decreases towards 0. Similarly (and in a symmetrical manner), the conditions are tight for $\theta \rightarrow 0$ when $b_N < |b_1|$, and can be relaxed increasingly as $|\theta|$ increases towards $\pi$. This follows from the proof, and it is demonstrated by numerical examples in \Cref{fig:ComplexEigRegions}. The blue region marks the domain where $r_{c^*}(b^c)\leq r^*$, that is, the results of \Cref{them:optimalC} hold so long as all the complex eigenvalues of $\mB$ lie within this region. The disk enclosed by the red circle is the subdomain covered by \Cref{them:Robustness}.

	
\subsection{Comparison of Nesterov's scheme to RI Chebyshev acceleration}
	
	A classical approach to accelerating convergence of stationary iteration schemes is by polynomial acceleration, whereby successive iterates are combined linearly with skillfully selected coefficients \cite{hageman2012applied}: 
	\begin{equation}
	\bar \vx_k = \sum_{n=0}^k\alpha_{k,n}\vx_n,\label{eq:PolynomialAcc}
	\end{equation}
	where $\vx_n = \mB\vx_{n-1}+\text{Constant}$ and $\sum_{n=0}^k\alpha_{k,n}=1$. Denoting $\bar \ve_k = \bar\vx_k-\vx^*$, we obtain
	$$
	\bar \ve_k = \left( \sum_{n=0}^k\alpha_{k,n} \mB^n\right)\ve_0 \, ,
	$$
	where $\ve_0=\vx_0-\vx^*$. The objective of minimizing the spectral radius $\rho\left( \sum_{n=0}^k\alpha_{k,n} \mB^n\right)$ (which yields the asymptotic convergence factor ACF), is achieved by using the well-known Chebyshev polynomial. Applied in recursive form, this yields a scheme of the following form:
	\begin{equation}
	\begin{array}{rcl}
	    \vx_1 &=& \gamma \left(\mB\vx_0+\text{Constant}\right)+(1-\gamma)\vx_0,\\
	    
	     \vx_{k+1}& =& {\beta}_{k+1}\left\{ \gamma \left(\mB\vx_k+\text{Constant}\right)+\left(1-\gamma\right)\vx_k\right\}+\left(1-\beta_{k+1}\right)\vx_{k-1}.
	\end{array}
	\label{eq:ChebyAccScheme}   
	\end{equation}
	The optimal values of $\gamma$ and $\beta_{k+1}$ depend on the eigenvalues of $\mB$, see details in \cite{hageman2012applied}. Often, such information is not readily available. This motivates us to consider a ``Restricted Information'' scenario, where we assume that we are given only the smallest and largest real eigenvalues of $\mB$, $b_1$ and $b_N$, as in the previous section. We refer to the scheme \eqref{eq:ChebyAccScheme} that is based solely on this information as Restricted Information (RI) Chebyshev acceleration. Clearly, RI Chebyshev acceleration is optimal in the case where all the eigenvalues of $\mB$ are real, because in that case we have all the required information. Indeed, in this case \eqref{eq:ChebyAccScheme} converges as fast as Preconditioned Conjugate Gradients (PCG) with preconditioner $\mB$, but requires less computation than PCG once $b_1$ and $b_N$ are known \cite{hageman2012applied}. In this case, Chebyshev acceleration is more efficient than Nesterov's scheme. However, it is interesting to note that we only need to know $b_N$ for applying Nesterov's scheme if $b_1\geq - \frac{1}{3}b_N$ (or, by a symmetric argument, we only need to know $b_1$ if $b_N\leq - \frac{1}{3}b_1$), so in these regimes less information is required---just the spectral radius of $\mB$ which is often easy to compute approximately. 
	
\begin{figure}[!htb]
		\centering
		\subfigure[$b_1=-0.3$ and $b_N=0.9$.]{ \includegraphics[scale=0.55]{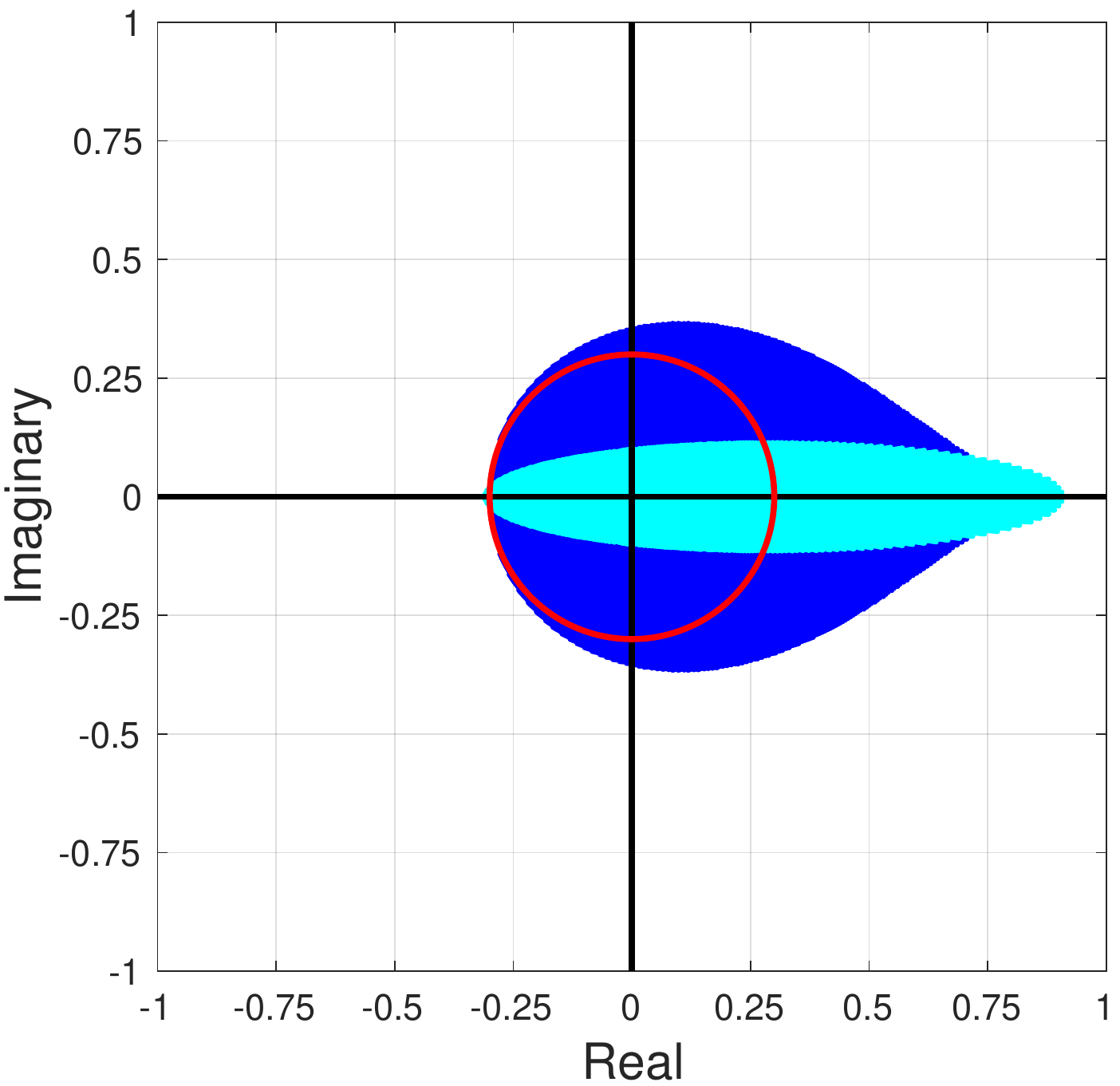}}
		 \hspace{0.2cm}
		\subfigure[$b_1=-0.5$ and $b_N=0.9$.]{ \includegraphics[scale=0.55]{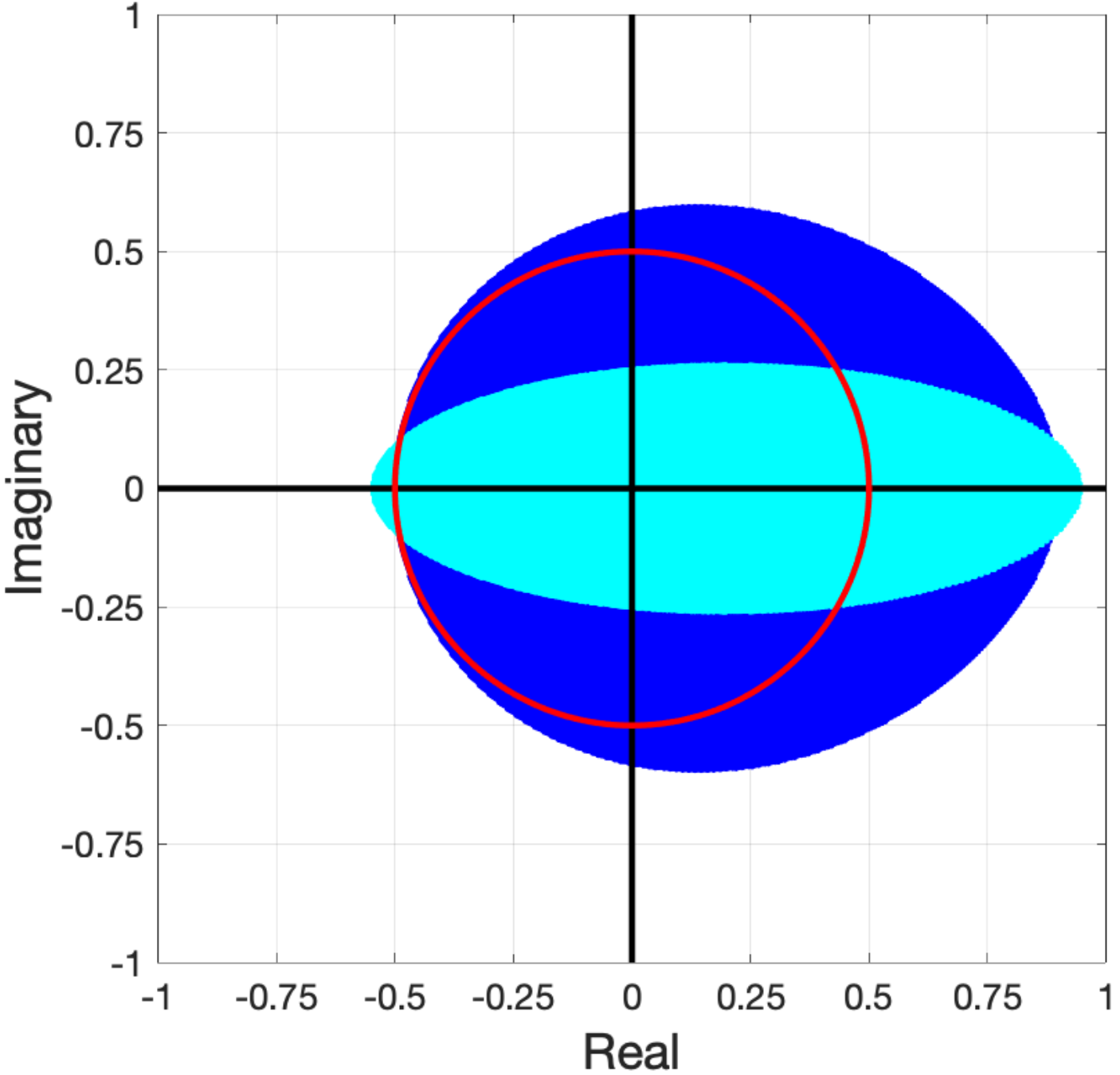}}
		\caption{Cyan: complex eigenvalues of $\mB$ for which RI Chebyshev acceleration yields a convergence factor smaller than or equal to $r^*$ of \Cref{them:optimalC}; Blue: complex eigenvalues of $\mB$ for which Nesterov's scheme yields a convergence factor smaller than or equal to $r^*$ of \Cref{them:optimalC}. The red circle is of radius $|b_1|$.}\label{fig:CompComplexEigRegions}
	\end{figure}

Next, we study and compare the performance of Nesterov's scheme and RI Chebyshev acceleration in cases where $\mB$ has complex eigenvalues\footnote{Optimal Chebyshev acceleration for $\mB$ with complex eigenvalues requires knowing an ellipse in the complex plane which contains all the eigenvalues, which may be hard to approximate in practice.}. As noted above, if $\mB$ has no complex eigenvalues then \eqref{eq:ChebyAccScheme} is always faster than Nesterov's scheme \cite{su2014differential}, as we also demonstrate later in our numerical tests. However, first we show in \Cref{fig:CompComplexEigRegions} that Nesterov's scheme can be faster than RI Chebyshev acceleration in a significant regime when $\mB$ does have complex eigenvalues. This figure compares Nesterov's scheme and RI Chebyshev for two cases of given $b_1$ and $b_N$ values. We use $r^*$, the ACF of \Cref{them:optimalC}, as a benchmark, and show the numerically computed range of eigenvalues in the complex plane for which RI Chebyshev (cyan) and Nesterov's scheme (blue) yield convergence factors which do not exceed $r^*$. Evidently, Nesterov's scheme remains robust for a significantly broader range of complex eigenvalues, and it converges faster than the RI Chebyshev acceleration scheme if there is at least one complex eigenvalue in the blue sub-domain. We demonstrate numerical examples of accelerated multigrid solvers where this behavior is relevant and yields an advantage to Nesterov's scheme. 
 
\begin{figure}[!htb]
		\centering
		\subfigure[$b_1=-0.3$ and $b_N=0.9$.]{ \includegraphics[scale=0.46]{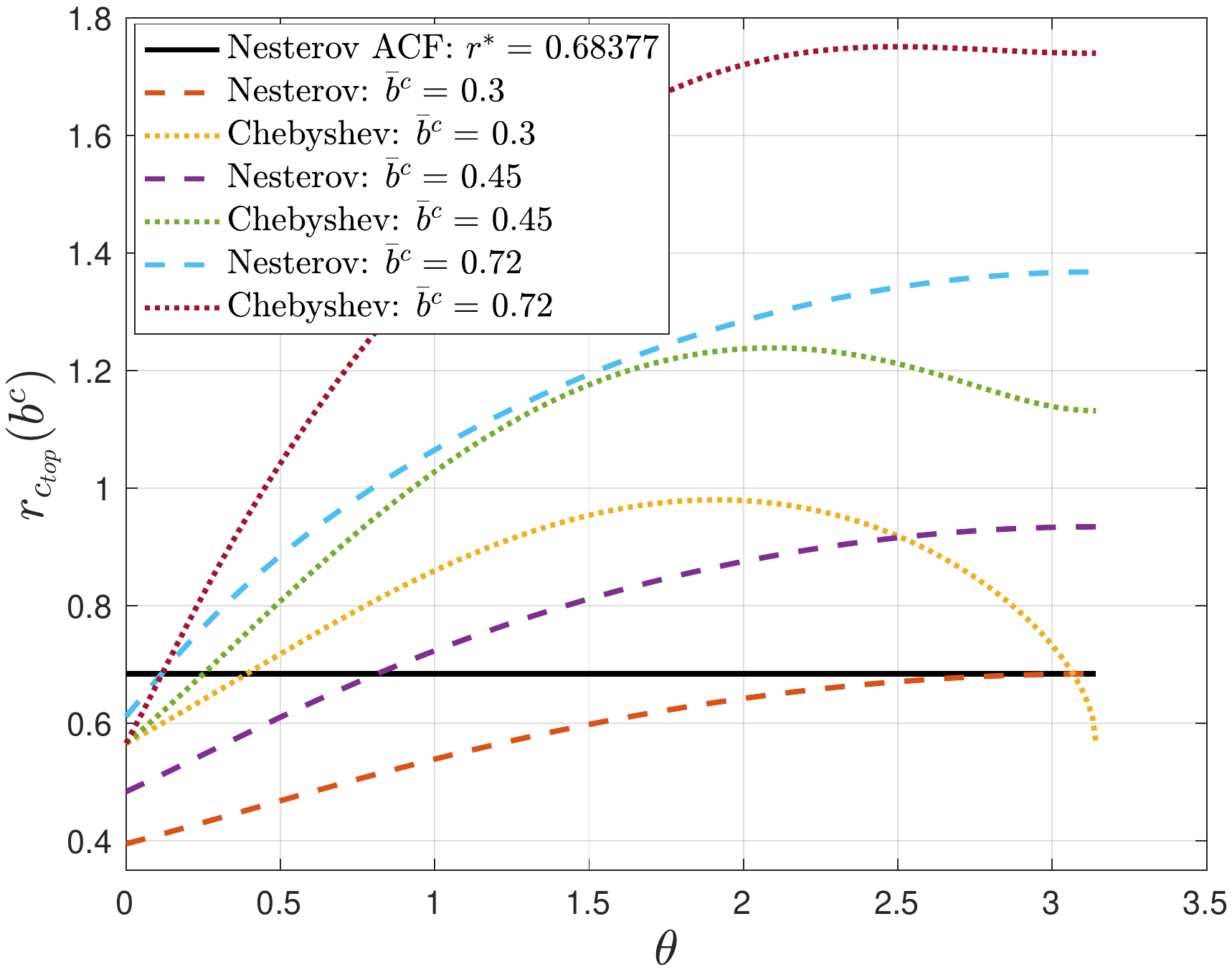}}
		\subfigure[$b_1=-0.5$ and $b_N=0.9$.]{ \includegraphics[scale=0.46]{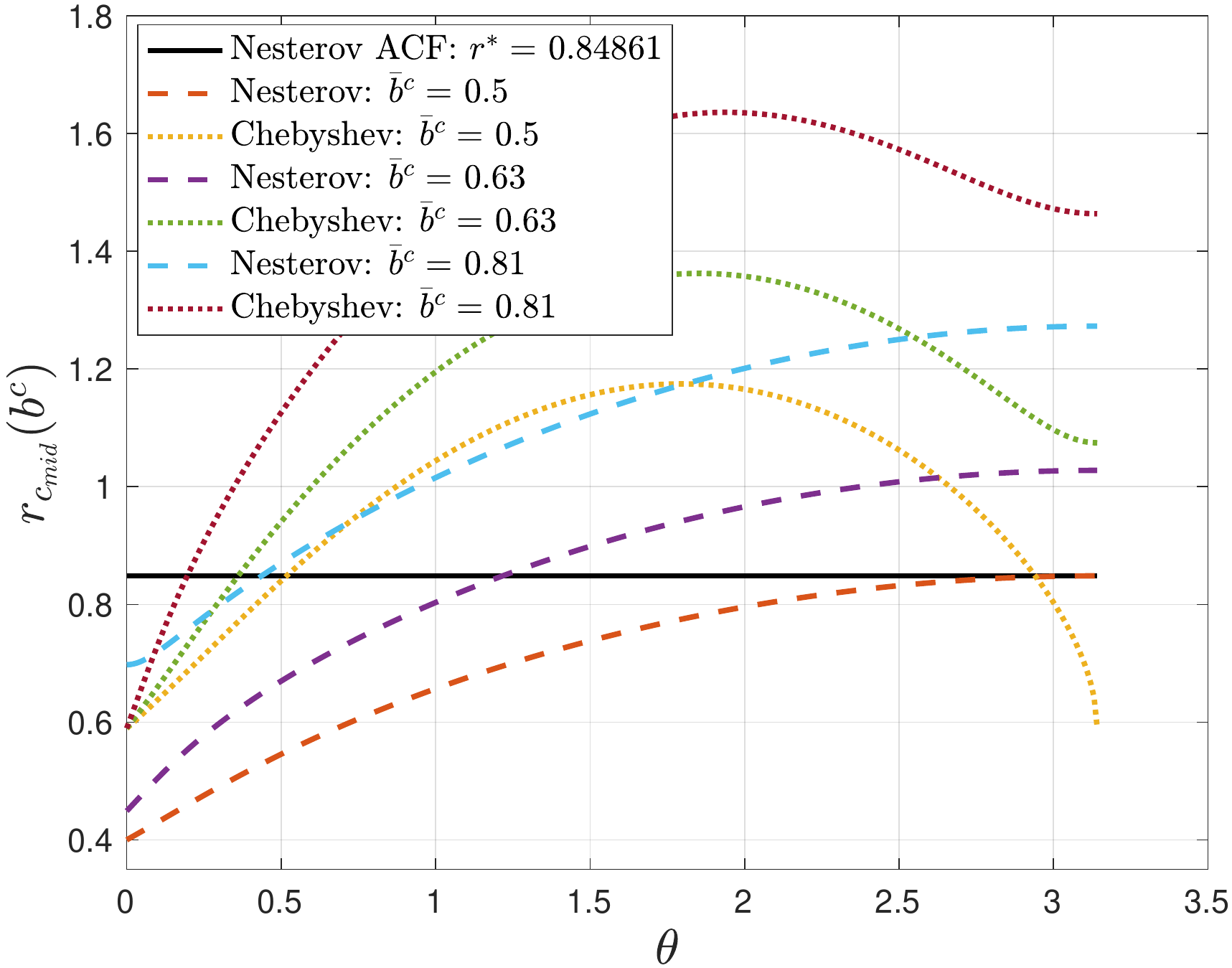}}
		\caption{ACF of particular complex eigenvalues $b^c$ with fixed $\bar b^c$ and varying $\theta\in[0,\pi]$ for Nesterov's scheme and RI Chebyshev acceleration. The black line represents $r^*$ of Nesterov's scheme. }\label{fig:CompComplexEigVStheta}
	\end{figure}

 Finally, we select several complex eigenvalues with a given $\bar b^c$ and argument varying from 0 to $\pi$, to show how the argument affects $r_{c^*}(b^c)$. The results are shown in \Cref{fig:CompComplexEigVStheta}. In general, we see that RI Chebyshev acceleration is adversely affected more strongly by change of $\theta$ than Nesterov's scheme, again demonstrating that the latter is more robust with respect to the existence of complex eigenvalues.   
	
\section{Numerical Tests}\label{sec:NumExps}

Nesterov's scheme is evidently easy to implement in practice for any stationary iteration method, since it only requires one additional step to combine the current iterate with the previous one. Furthermore, the additional computation is negligible and so acceleration is obtained almost for free. The only significant cost is the memory, since we need to store one previous iterate. In particular, compared with common acceleration techniques such as Krylov subspace methods \cite{saad2003iterative}, the cost of \eqref{eq:NesterovLinear:General} is smaller. The drawback of course is the requirement to know $b_1$ and $b_N$, since they are needed for computing the optimal parameter  $c^*$. In practice, as noted earlier, we may only need to know the spectral radius of $\mB$. For example, if $b_1\geq 0$ (e.g., if we iterate with $\mB$ twice before successive Nesterov steps), then we only need to evaluate $b_N$, which can be done approximately by the power method \cite{golub2000eigenvalue}, or by running the unaccelerated iteration for several steps. Since in many applications we need to solve \eqref{eq:linearPro} many times with different $\vf$, we argue that the amount of computation required to approximate $b_1$ and $b_N$ is often acceptable \cite{saad2003iterative}.

In this section we focus on accelerating multigrid V-cycle iterations \cite{brandt1977multi,briggs2000multigrid,yavneh2006multigrid} for elliptic boundary value problems. In some applications, the so-called smoothing factor, which is obtained by Fourier smoothing analysis, may provide a cheap yet sufficiently accurate approximation to $b_1$ and $b_N$  \cite{trottenberg2000multigrid,wienands2004practical,yavneh1995multigrid,yavneh1996red,yavneh1998multigrid}. 

We compare acceleration schemes for a multigrid solver for the two-dimensional diffusion equation on the unit square,
    \begin{equation}
	    -\nabla (\sigma(x,y) \nabla u(x,y)) = f(x,y),\label{eq:diffusionproblems}
	\end{equation} 
	discretized on a $1024\times1024$ uniform grid, yielding a linear system 
	\begin{equation}
	    \mA \vu = \vf \, , \label{eq:DiscreteDiffusion}
	\end{equation} 
	with $\mA = -\nabla^h (\vsigma \nabla^h \vu)$ a second-order finite-difference or bilinear finite element  discretization of the diffusion equation with Dirichlet boundary conditions. In the first set of tests, the elements of the diffusion coefficient vector are all equal to 1, yielding a discretized Poisson equation, while in the second and third sets the elements are sampled from a log-normal distribution and from a uniform distribution in $(0,1)$, respectively, following \cite{greenfeld2019learning}. All the tests in this paper are implemented on a laptop with $2.3$GHz Intel Core i$9$.
	
	For the Poisson problem we employ the standard five-point finite difference discretization, damped Jacobi or Red-Black Gauss-Seidel relaxation, full-weighted residual transfers and bilinear interpolation, and the coarse-grid operators are defined by rediscretization on all coarse grids with the five-point discretization stencil \cite{trottenberg2000multigrid}. We run standard V$(\nu_1,\nu_2)$ cycles with $\nu_1 = 1$ pre-relaxation sweep and $\nu_2 = 0$ or $1$ post-relaxation sweep. We compare acceleration by Nesterov's scheme to Preconditioned Conjugate Gradients (PCG) with V$(\nu_1,\nu_2)$ as the preconditioner, denoted PCG-V$(\nu_1,\nu_2)$ and to RI Chebyshev acceleration of the V$(\nu_1,\nu_2)$ cycles, denoted Chebyshev-V$(\nu_1,\nu_2)$.

    Let $\vr_k = \vf - \mA \vu_k$ denote the residual vector at the end of the $k$th iteration. Then the convergence factor CF at the $k$th iteration is given by the ratio of the successive residual norms, $||\vr_k||_2 / ||\vr_{k-1}||_2$. We estimate the ACF by running sufficiently many iterations such that the CF no longer changes significantly. 
	
\begin{figure}[!htb]
		\centering
	    
	    \includegraphics[scale=0.5]{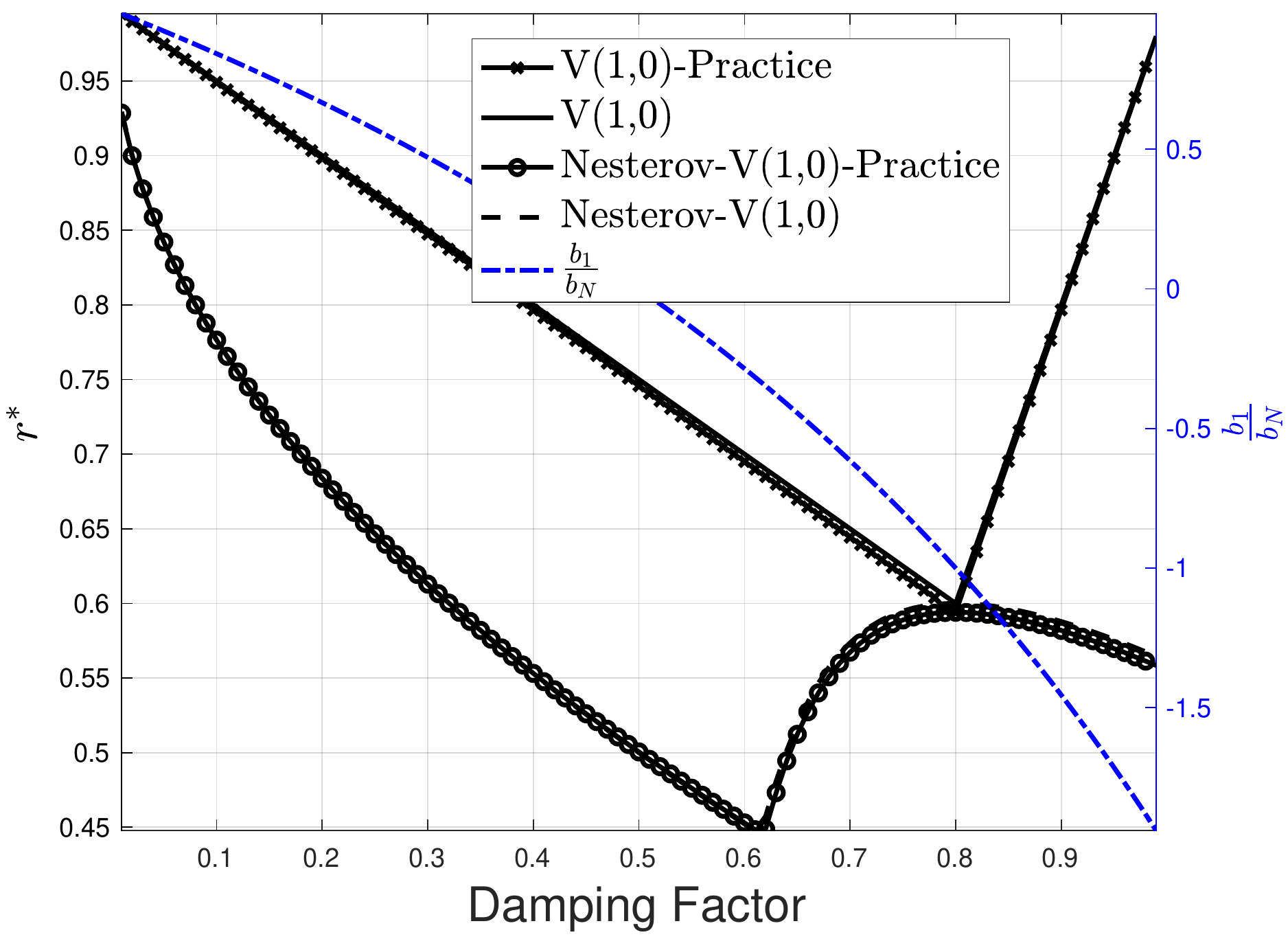}
		\caption{The ACF achieved by V$(1,0)$ cycles, with and without Nesterov acceleration, as a function of the Jacobi relaxation damping factor. The extreme eigenvalues used for determining $c^*$, $b_1$ and $b_N$, are estimated by Fourier smoothing analysis, and the ratio $\frac{b_1}{b_N}$ is also shown. ``Practice'' refers to results achieved in practice by numerical computations on a $256\times 256$ grid, terminated when the residual norm is smaller than $10^{-8}$. The ACF is then estimated by the geometric mean of the last $5$ iterations.}
		\label{fig:ConvRateRegionmark}
	\end{figure}

\begin{figure}[!htb]
		\centering
		\subfigure[]{\includegraphics[scale=0.3]{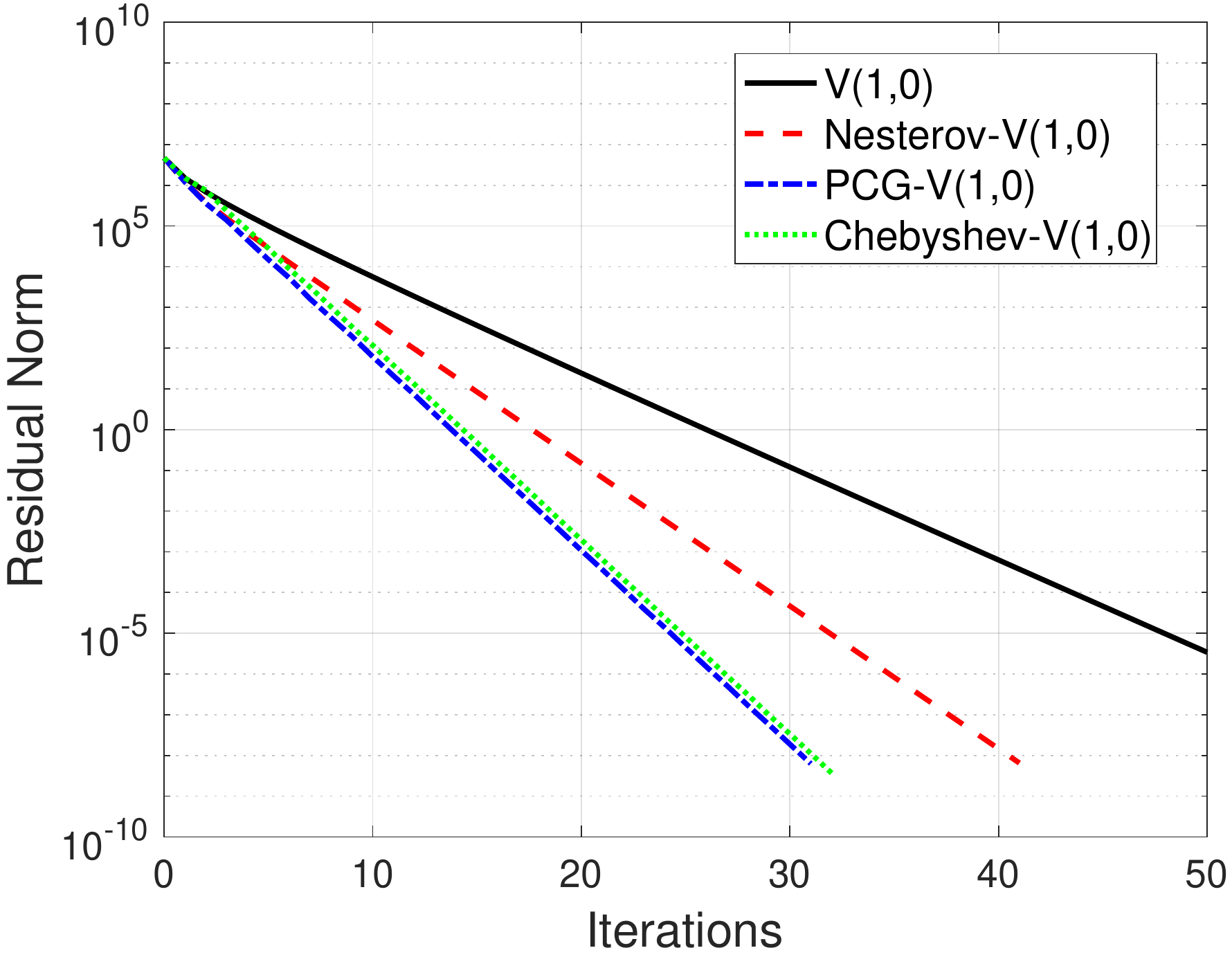}\label{fig:PoissonExampJacobiResi10}}
		\subfigure[]{\includegraphics[scale=0.3]{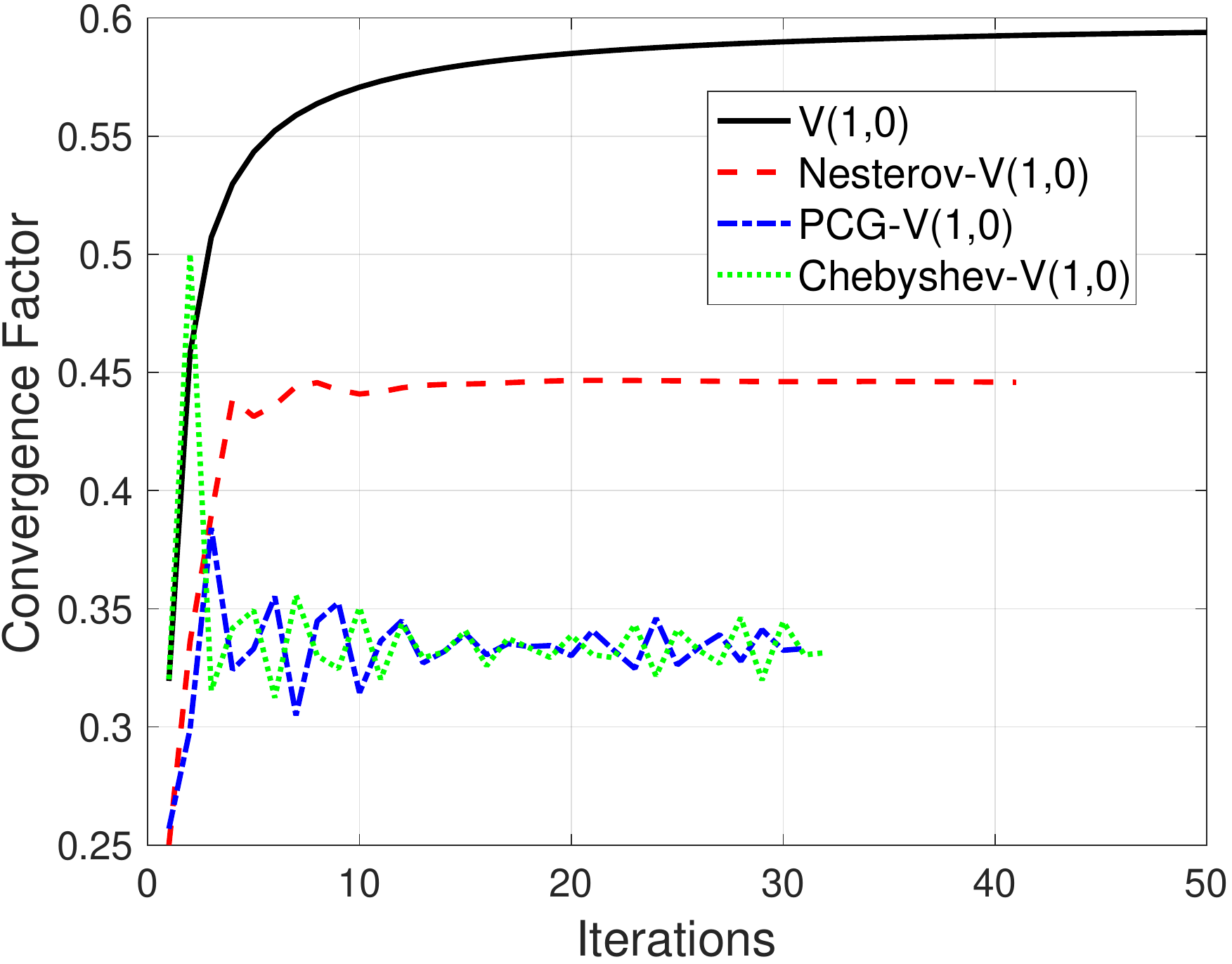}\label{fig:PoissonExampJacobiConv10}}
		\subfigure[]{\includegraphics[scale=0.3]{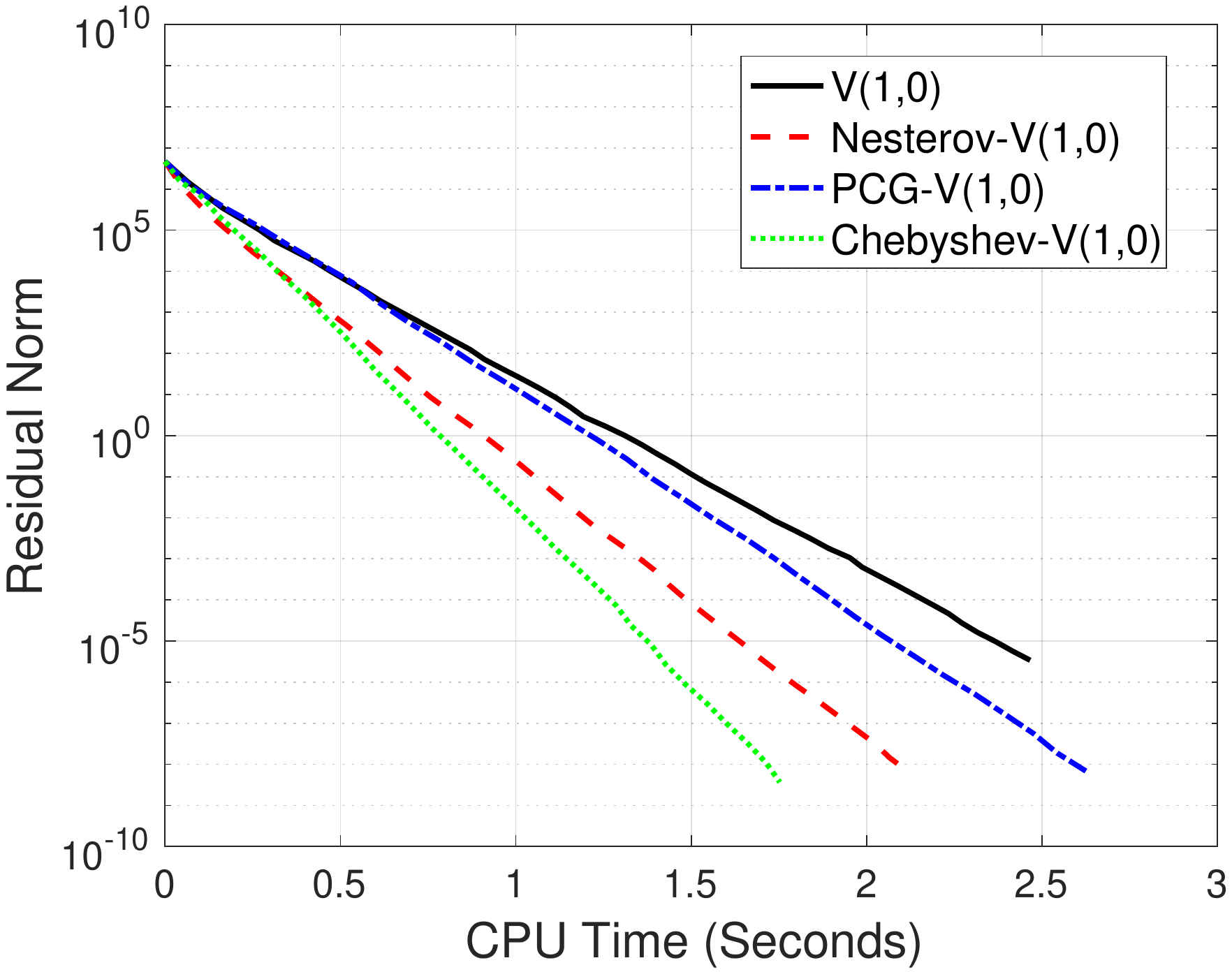}\label{fig:PoissonExampJacobiCPU10}}
		
		\subfigure[]{\includegraphics[scale=0.3]{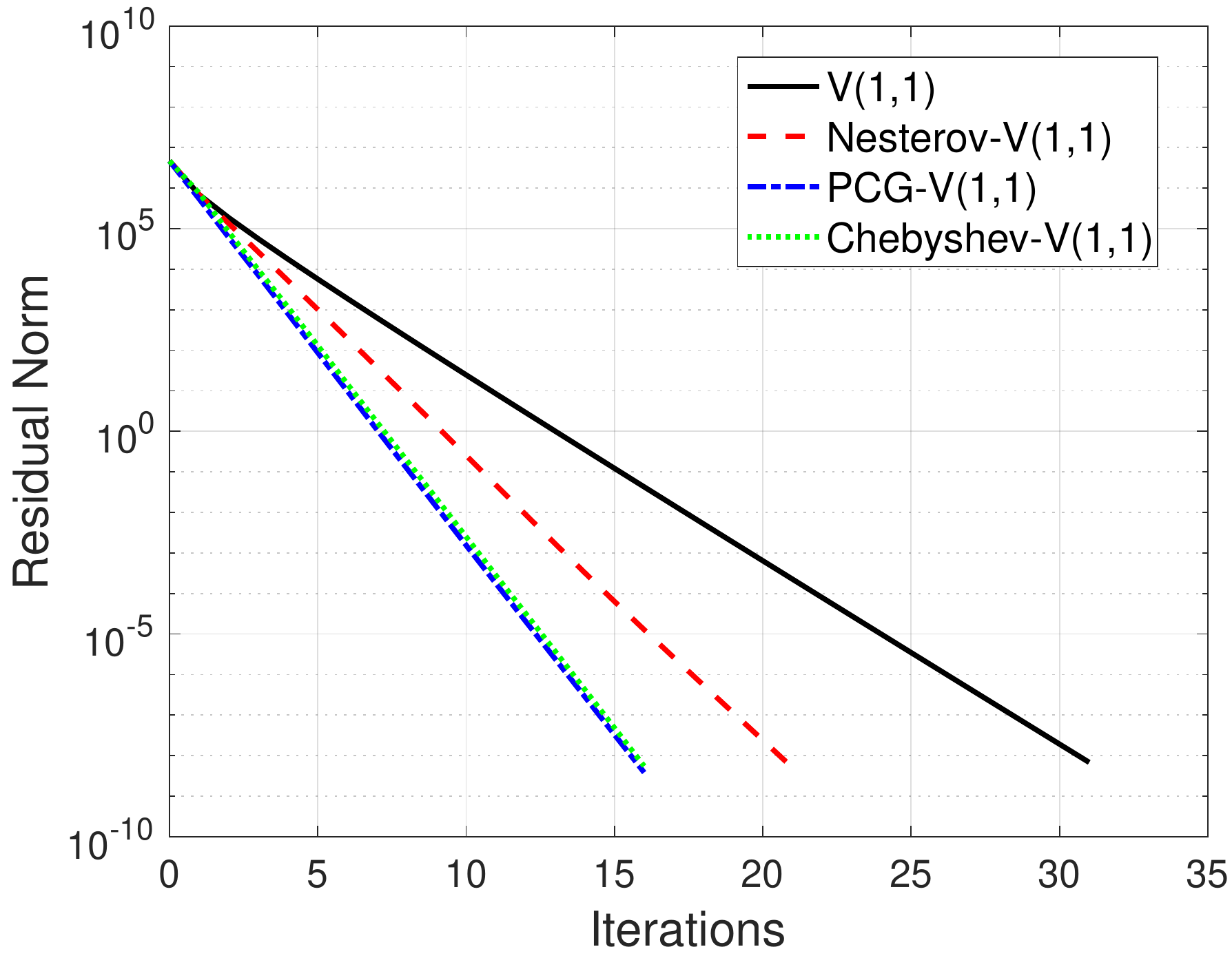}\label{fig:PoissonExampJacobiResi11}}
		\subfigure[]{\includegraphics[scale=0.3]{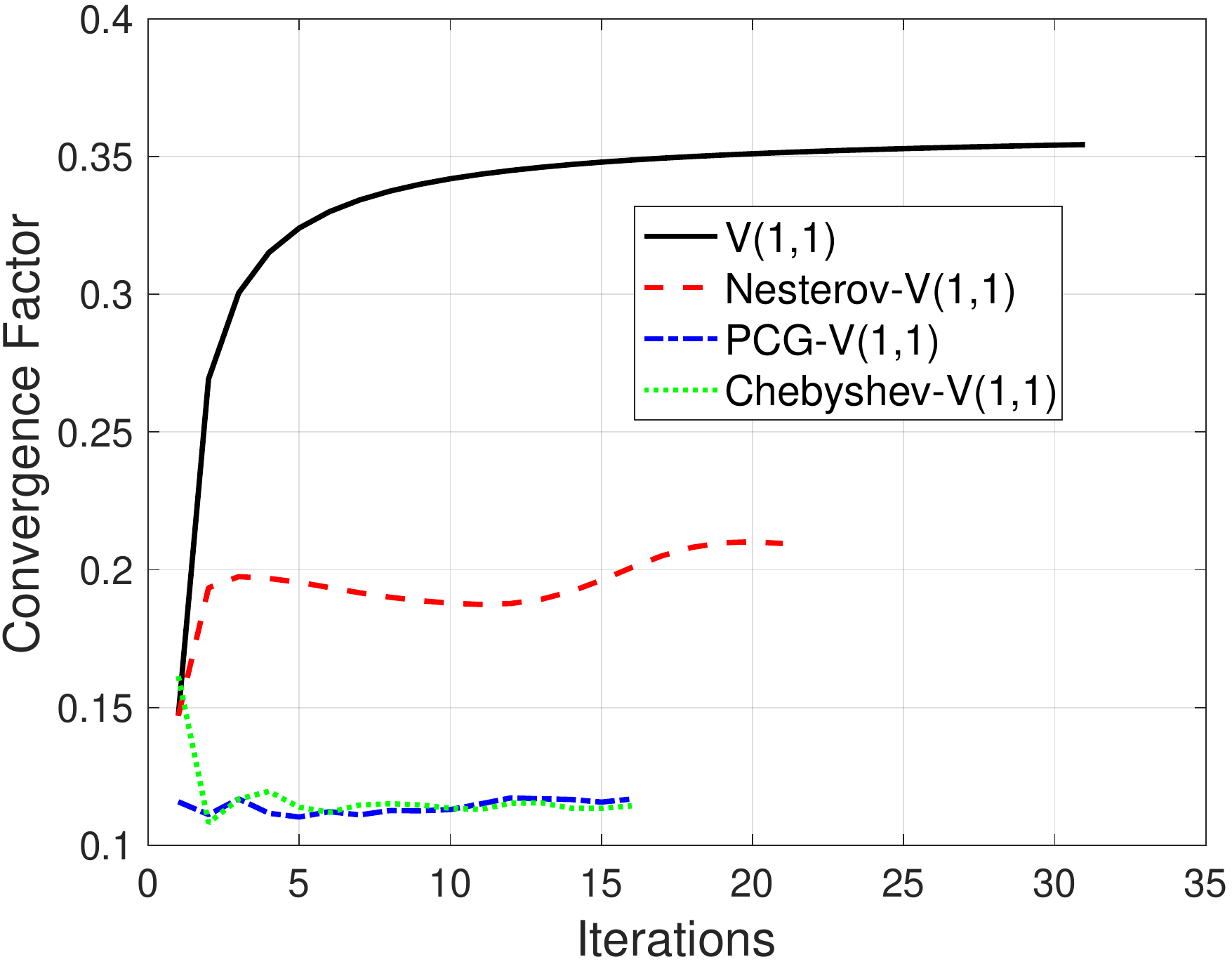}\label{fig:PoissonExampJacobiConv11}}
		\subfigure[]{\includegraphics[scale=0.3]{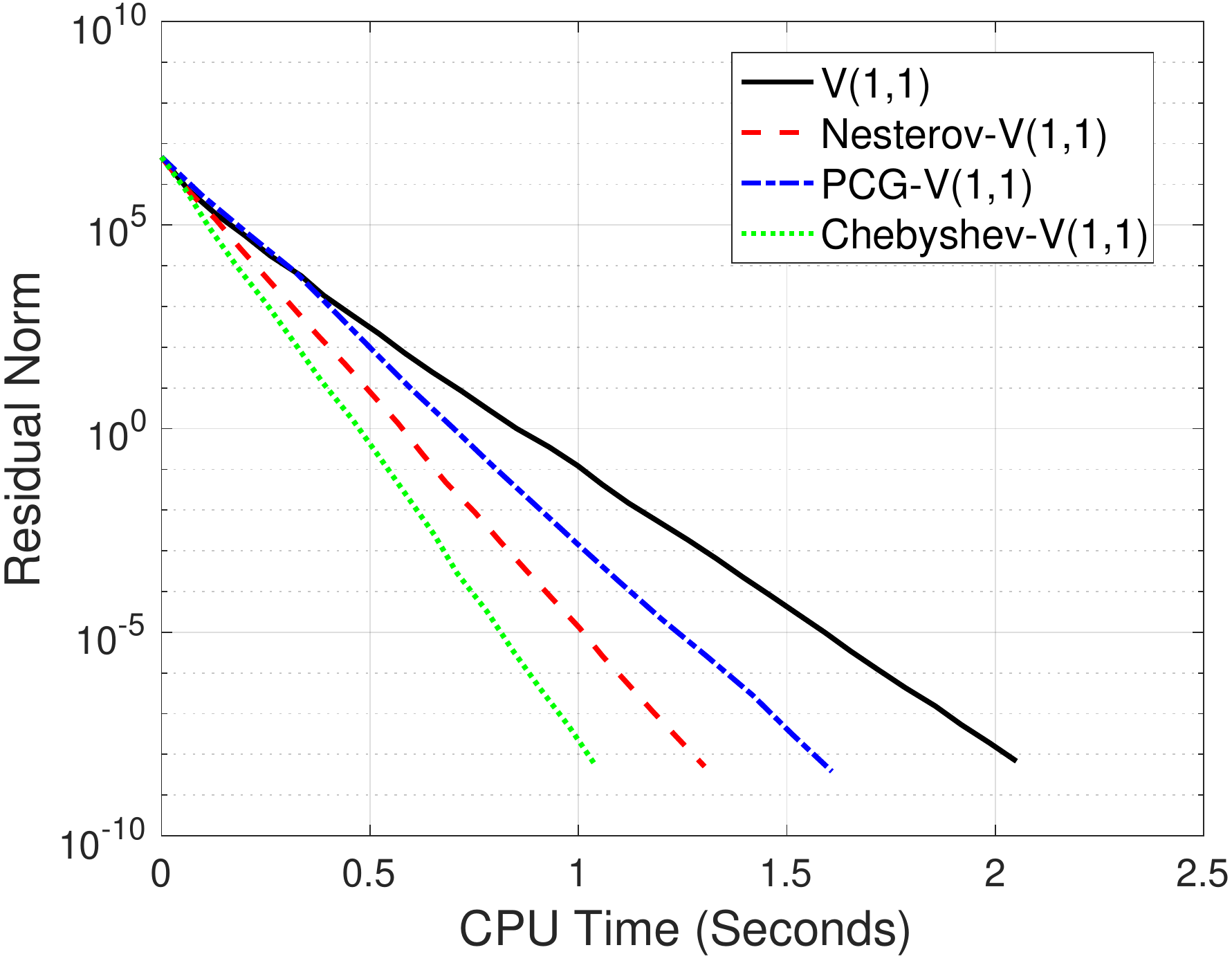}\label{fig:PoissonExampJacobiCPU11}}
		\caption{Comparison of acceleration methods for the Poisson problem. Optimally damped Jacobi relaxation is used in all the tests except for Nesterov-V$(1,0)$, which uses a damping coefficient of $\frac{8}{13}$. First row: Accelerated V$(1,0)$ cycles. Second row: Accelerated V$(1,1)$ cycles.}
		\label{fig:PoissonExampJacobi}
	\end{figure}

	In our first test we employ V(1,0) cycles for the Poisson problem, using Jacobi relaxation with the theoretical optimal damping factor 0.8 (obtained by Fourier smoothing analysis), yielding a smoothing factor (and correspondingly, an ACF) of 0.6. For this choice, Nesterov's scheme cannot provide acceleration, because $b_N = -b_1 = 0.6$. However, it is possible to improve the convergence factor by choosing a damping factor that is not optimal for stand-alone multigrid. We do this by increasing both $b_N$ and $b_1$, and of course applying Nesterov acceleration. To demonstrate the potential gain, we show in \Cref{fig:ConvRateRegionmark} the ACF of V$(1,0)$ and Nesterov-accelerated V$(1,0)$ cycles, for varying damping factors of the Jacobi relaxation. We find that the optimal choice is to reduce the damping factor to $\frac{8}{13}$, yielding $b_1=-\frac{1}{3}b_N$. Evidently, these theoretical results are matched well by the numerical results achieved in practice. The effect is also seen in the first row of \Cref{fig:PoissonExampJacobi}. We find that optimizing the Jacobi damping factor for Nesterov acceleration yields a non-negligible reduction in the ACF (from 0.6 to 0.45) and in the run-time. Moreover, although PCG-V$(1,0)$ is faster than Nesterov-V$(1,0)$ in terms of iteration count, as expected, Nesterov acceleration yields a shorter run-time. As expected, the winner is Chebyshev-V$(1,0)$, which is as fast as PCG-V$(1,0)$ in terms of iteration count, but faster than all its competitors in terms of CPU time. Rather similar conclusions are obtained for V$(1,1)$ cycles, as seen in the second row of \Cref{fig:PoissonExampJacobi}.
	
	
\begin{figure}[!htb]
	    \centering
	    \subfigure[]{\includegraphics[scale=0.3]{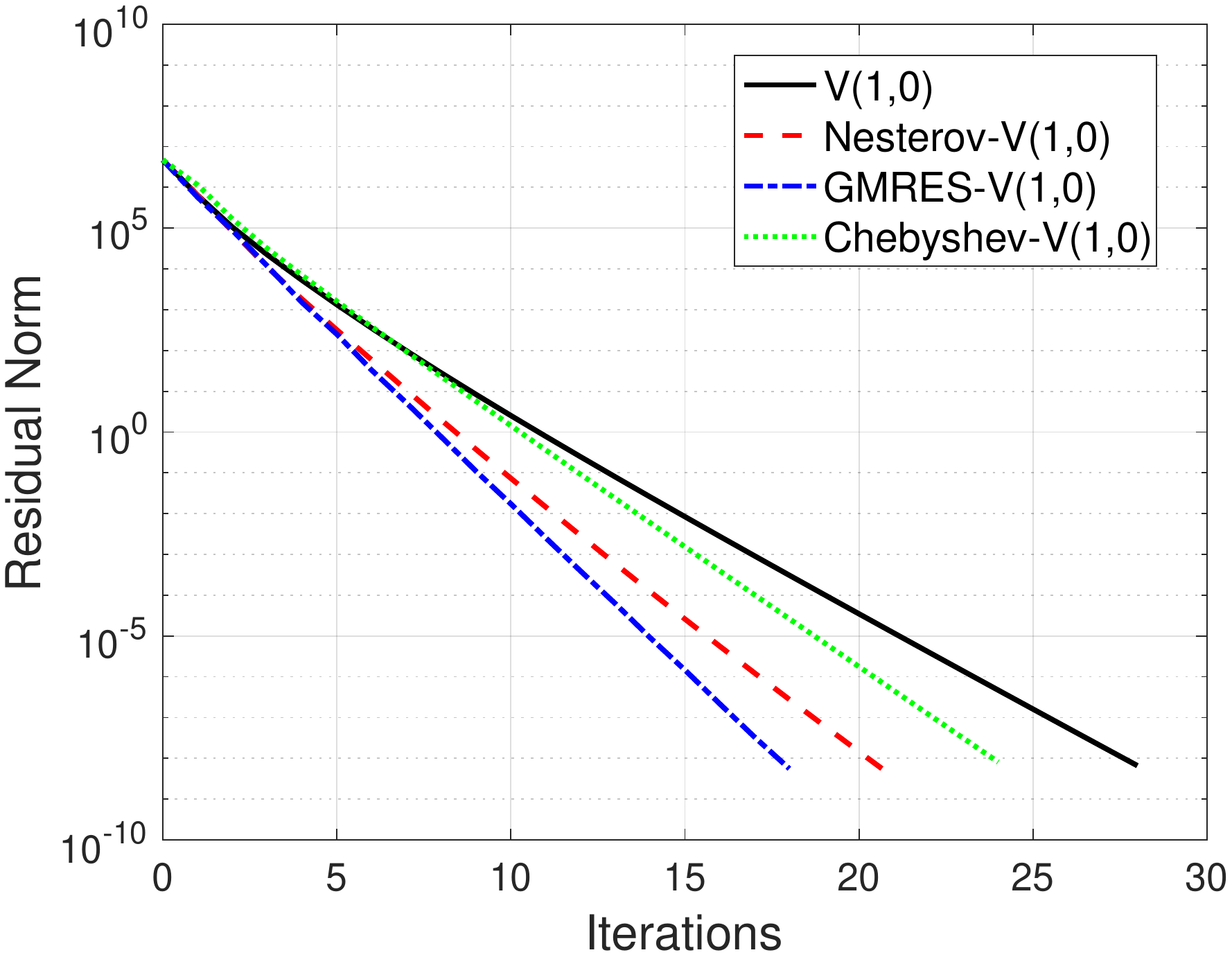}\label{fig:PoissonExampRBResi10}}
	    \subfigure[]{\includegraphics[scale=0.3]{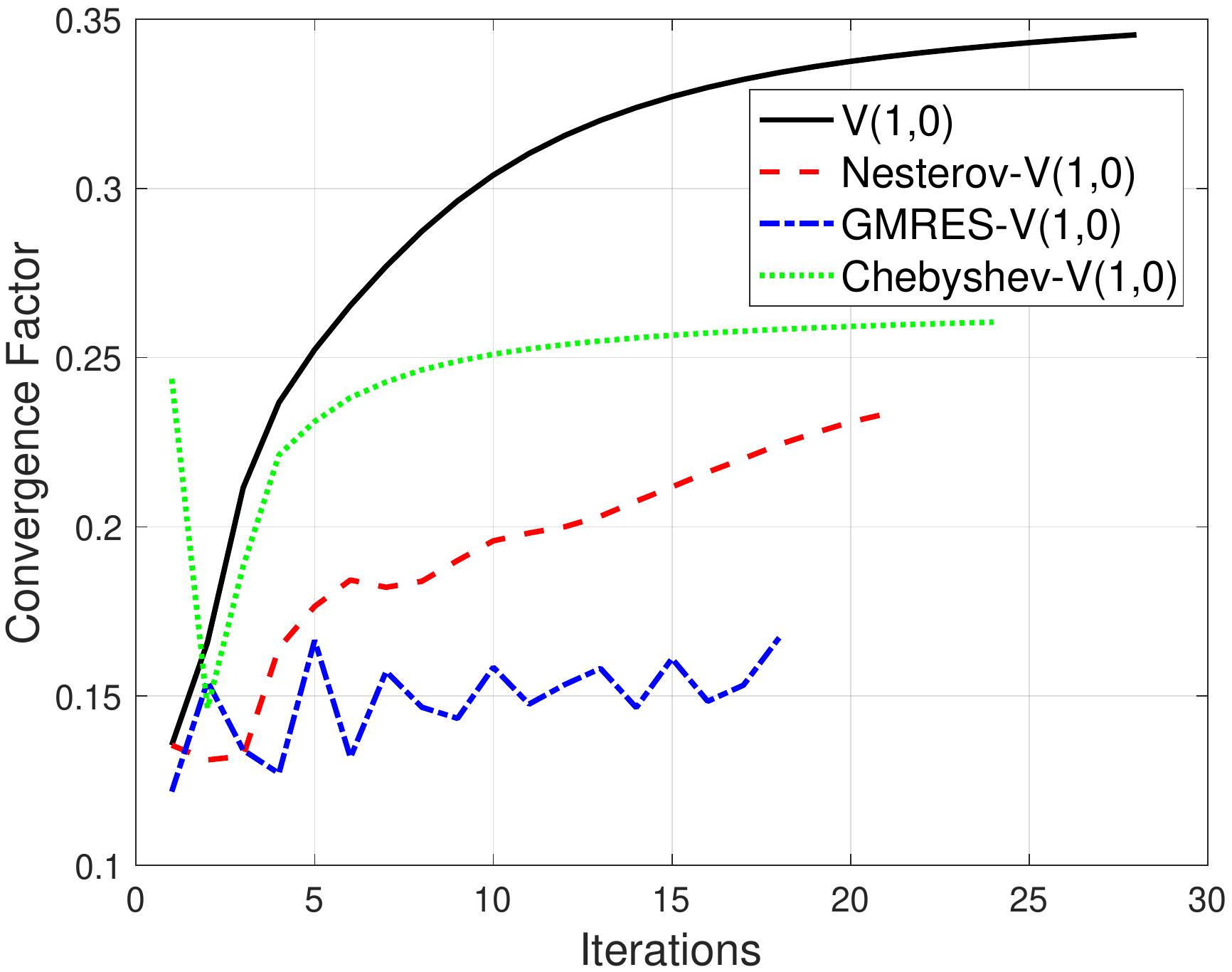}\label{fig:PoissonExampRBConv10}}
	    \subfigure[]{\includegraphics[scale=0.3]{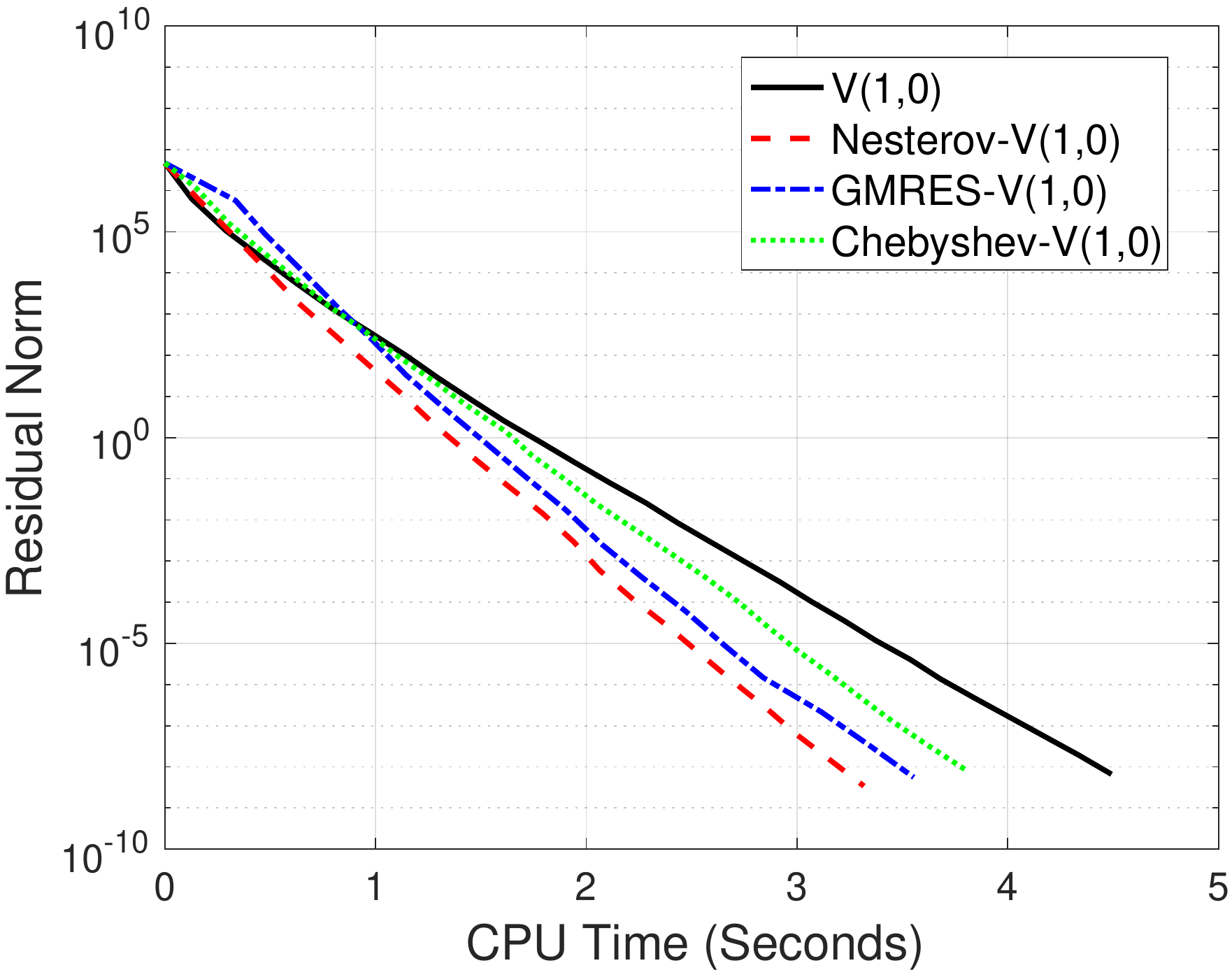}\label{fig:PoissonExampRBCPU10}}
	     
	    \subfigure[]{\includegraphics[scale=0.3]{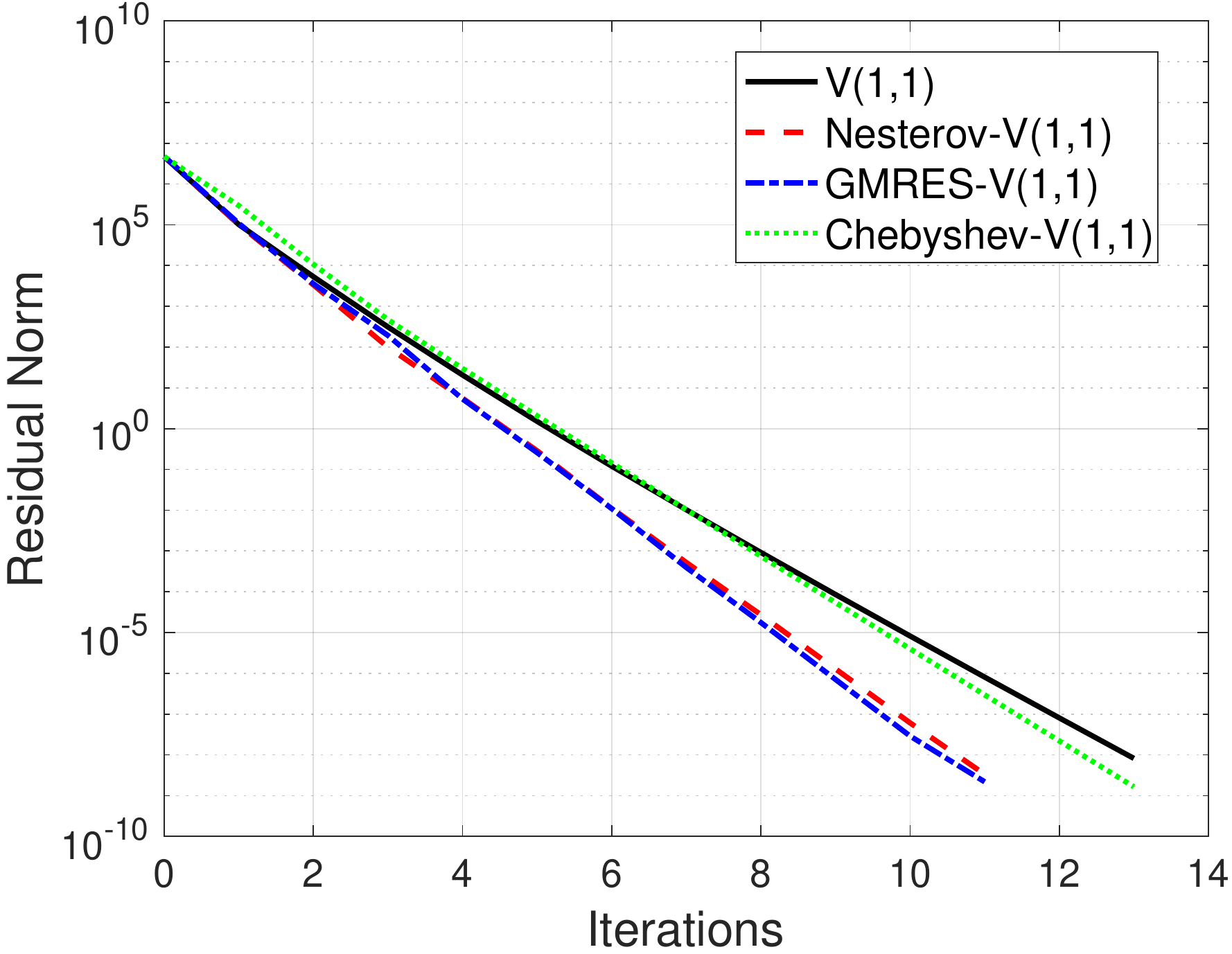}\label{fig:PoissonExampRBResi11}}
	    \subfigure[]{\includegraphics[scale=0.3]{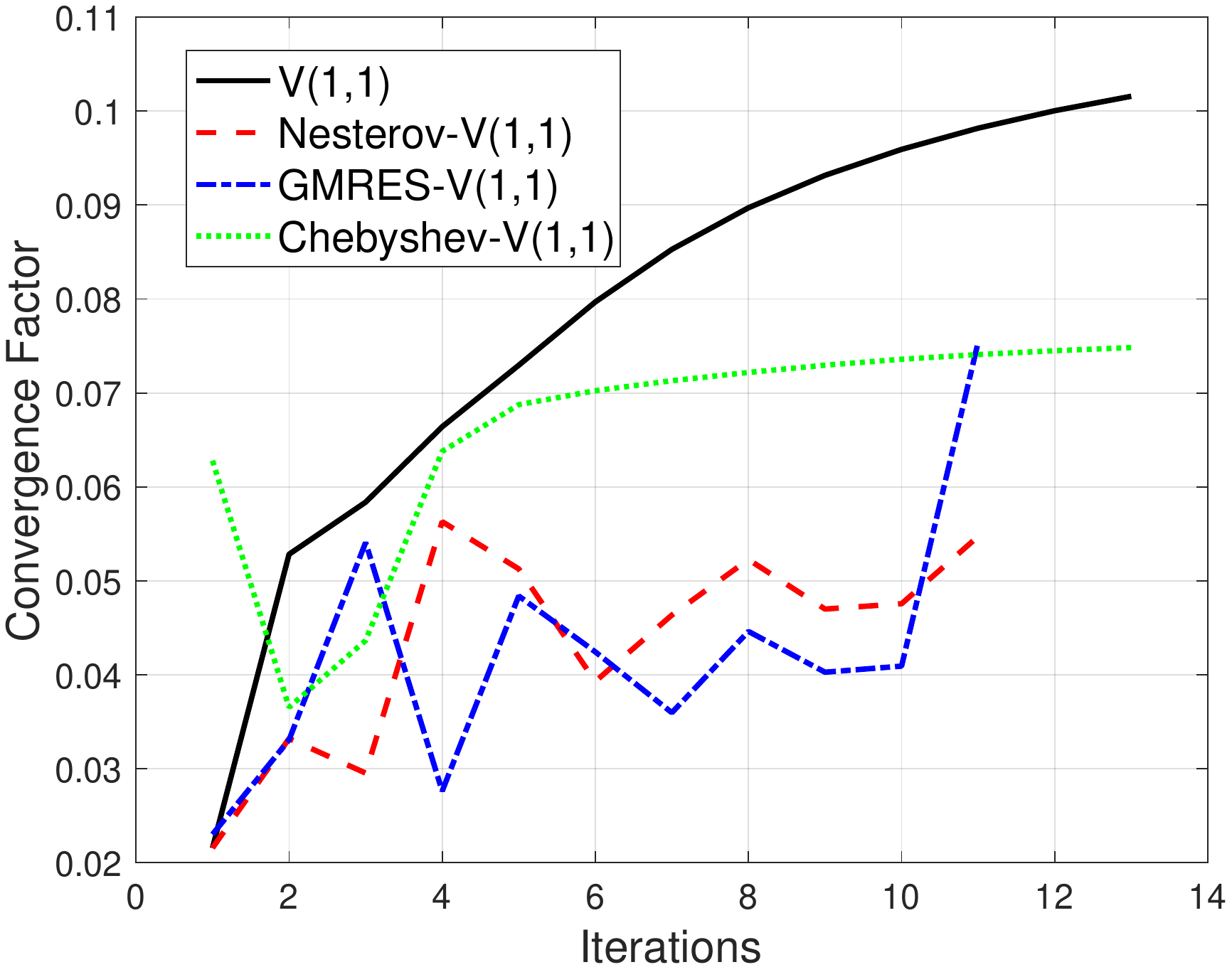}\label{fig:PoissonExampRBConv11}}
	    \subfigure[]{\includegraphics[scale=0.3]{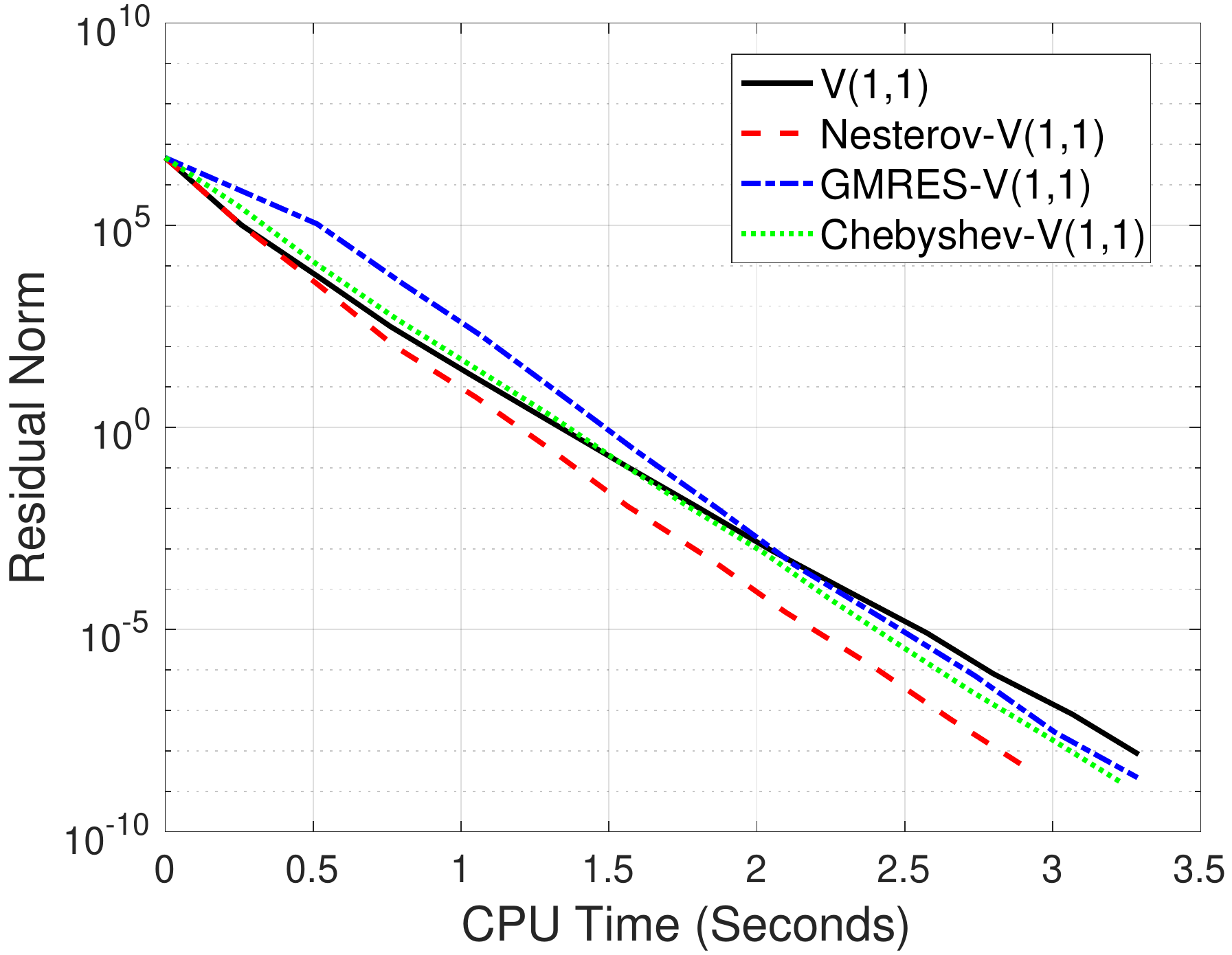}\label{fig:PoissonExampRBCPU11}}
	    \caption{Comparison of acceleration methods for the Poisson problem with Red-Black relaxation. First row: Accelerated V$(1,0)$ cycles. Second row: Accelerated V$(1,1)$ cycles.}
	    \label{fig:PoissonExampRB}
	\end{figure}

	We next test accelerated multigrid cycles with Red-Black (RB) relaxation, which costs about the same as Jacobi and provides better smoothing, hence faster convergence. In this case, some of the eigenvalues of $\mB$ are complex, and yet we select $c^*$ based only on $b_1$ and $b_N$ using \Cref{them:Robustness}. Again, we test both V$(1,0)$ and V$(1,1)$ cycles with acceleration. The results are shown in \Cref{fig:PoissonExampRB}. Note that we use GMRES without restart \cite{saad2003iterative}, {\color{black} which is faster than restarted GMRES}, instead of PCG in this experiment, because V$(1,0)$ and V$(1,1)$ are not symmetric. {\color{black}The alternative of using symmetric Gauss-Seidel for the relaxation instead of RB retains the symmetry, allowing the use of PCG. However, in our numerical tests with $\nu_1=\nu_2=1$ we found that using symmetric GS yields an ACF that is similar to that of RB but requires twice the number of operations per cycle. Moreover, RB has an advantage in parallel computation.} In \Cref{fig:PoissonExampRB}, we see that Nesterov's scheme is faster than RI-Chebyshev acceleration, both in terms of iterations and CPU time. Moreover, we observe that GMRES with V$(1,0)$ or V$(1,1)$ as the preconditioner is fastest in terms of iterations. However, Nesterov's scheme is fastest in terms of CPU time. Furthermore, we note that GMRES needs much more memory than Nesterov's scheme, and its implementation is not as trivial.
	
\begin{figure}[!htb]
		\centering
		\subfigure[]{\includegraphics[scale=0.3]{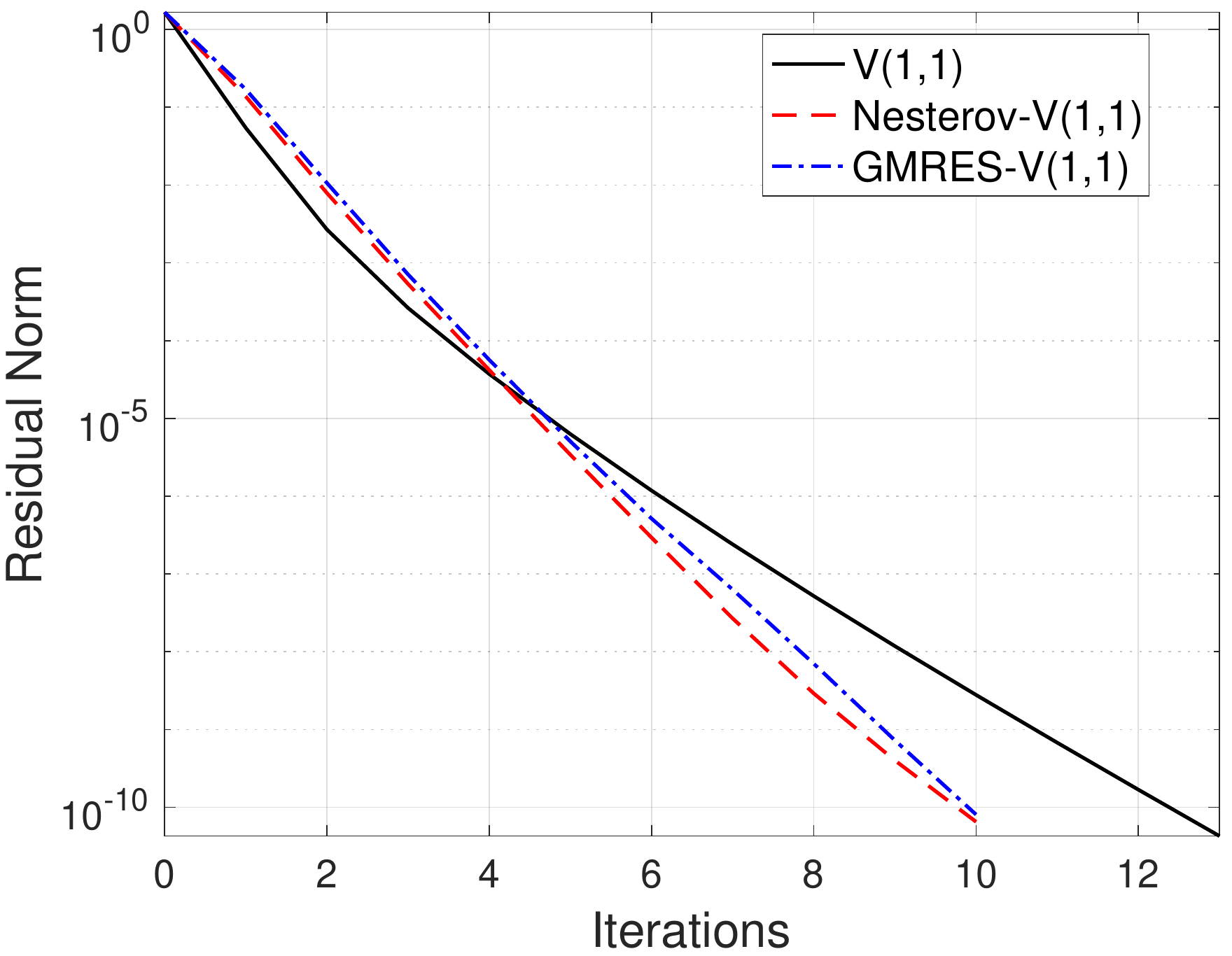}\label{fig:DiffusionExample_lognormIter}}
		\subfigure[]{\includegraphics[scale=0.3]{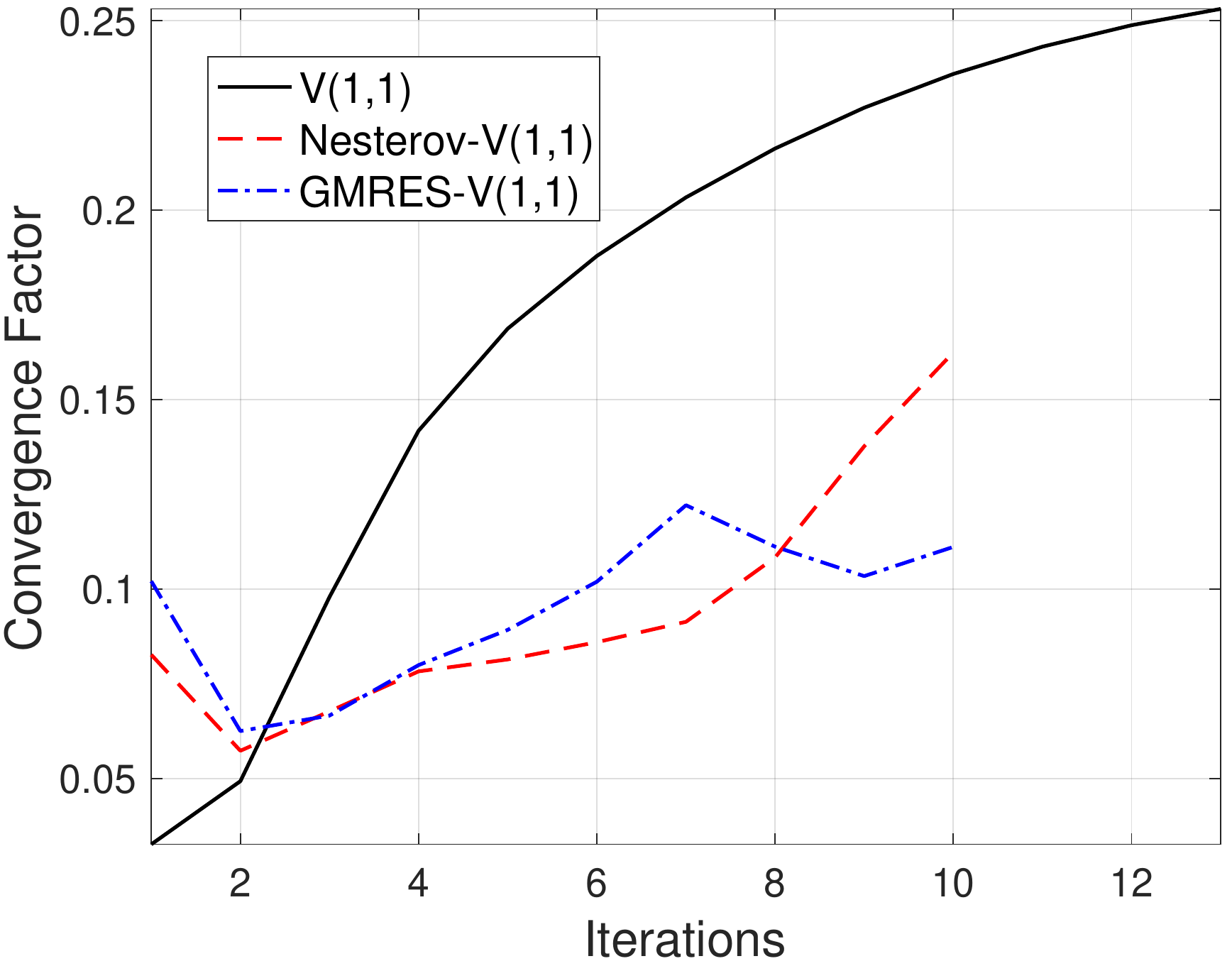}\label{fig:DiffusionExample_lognormConv}}
		\subfigure[]{\includegraphics[scale=0.3]{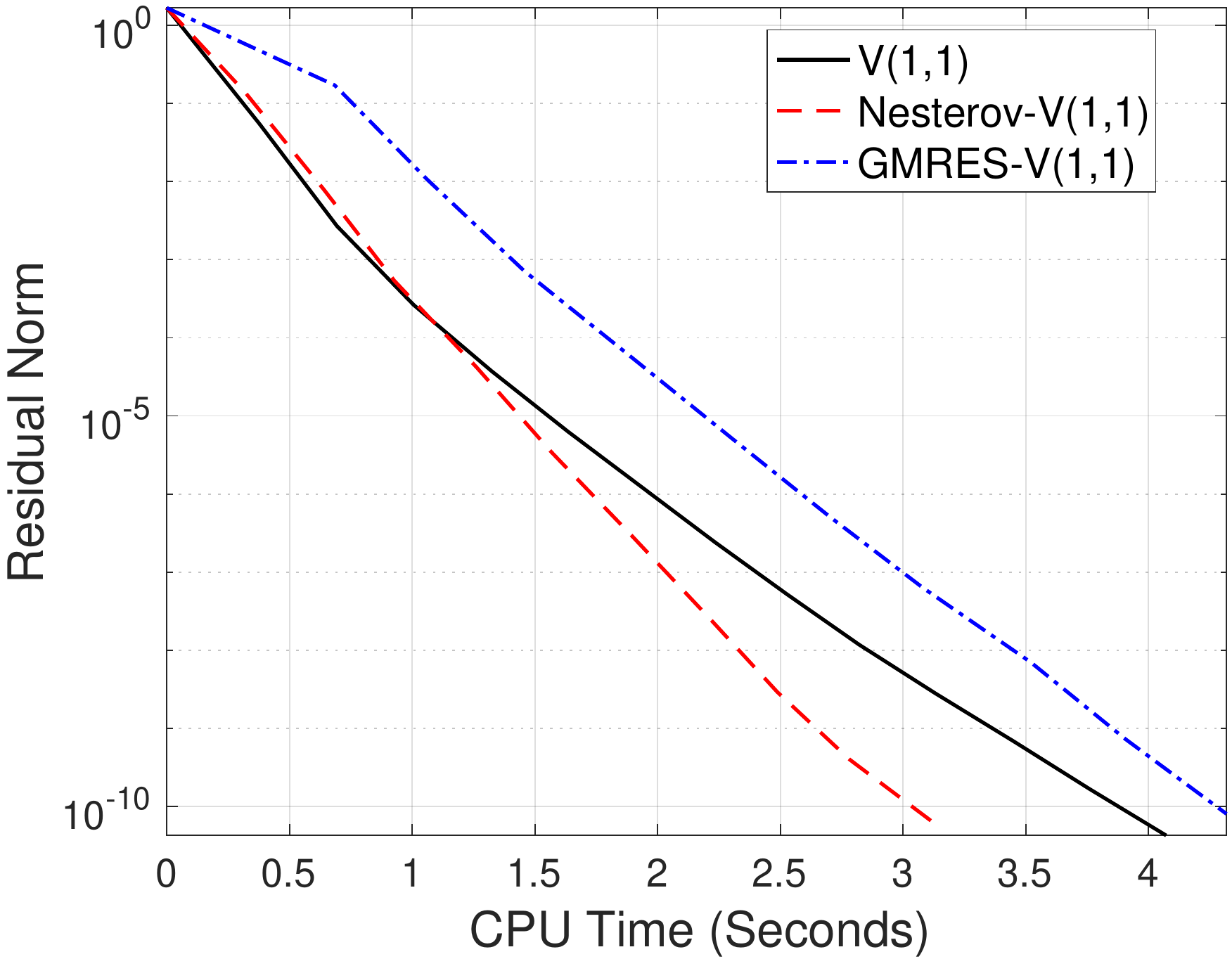}\label{fig:DiffusionExample_lognormCPU}}
		
		\subfigure[]{\includegraphics[scale=0.3]{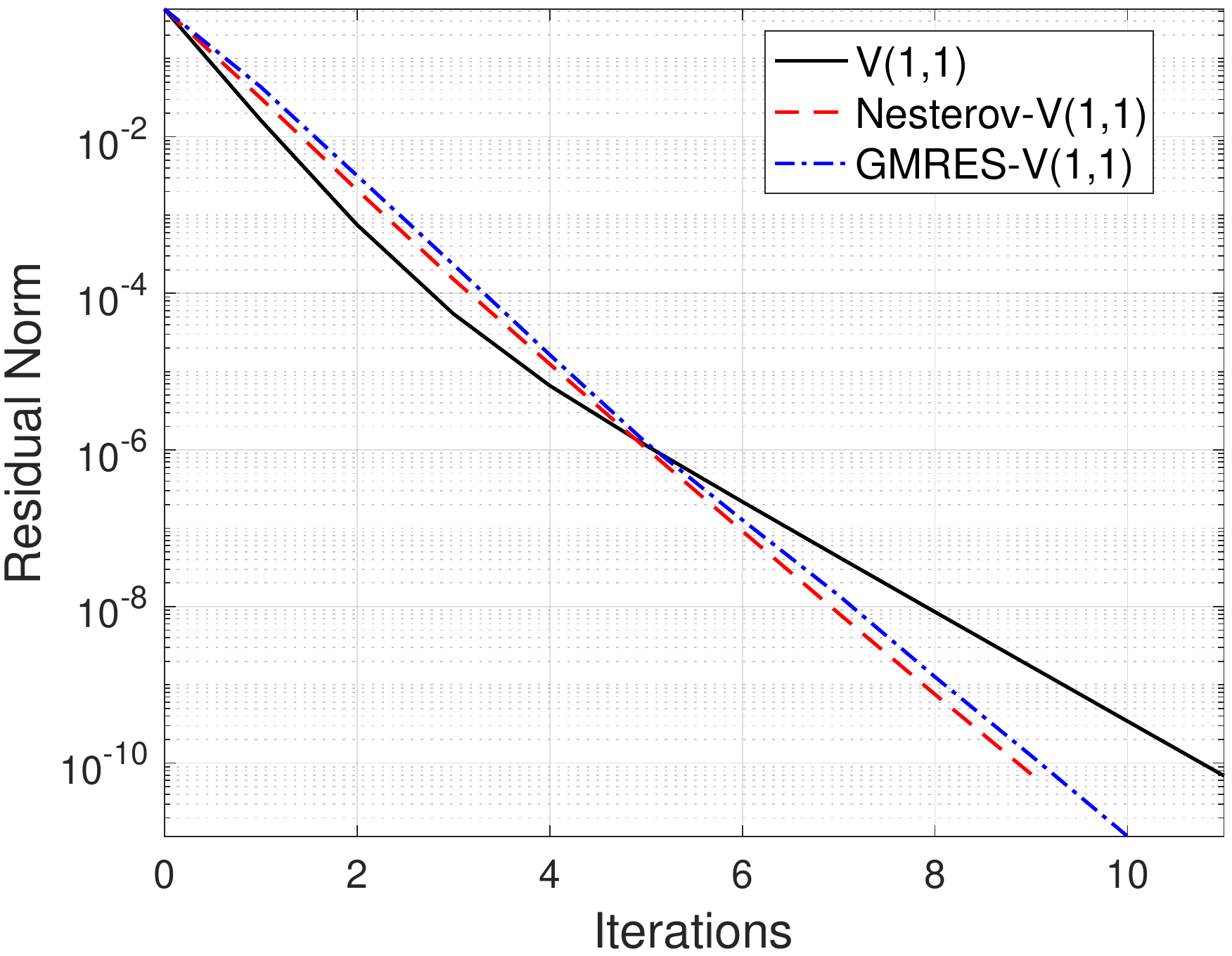}\label{fig:DiffusionExample_uniformIter}}
		\subfigure[]{\includegraphics[scale=0.3]{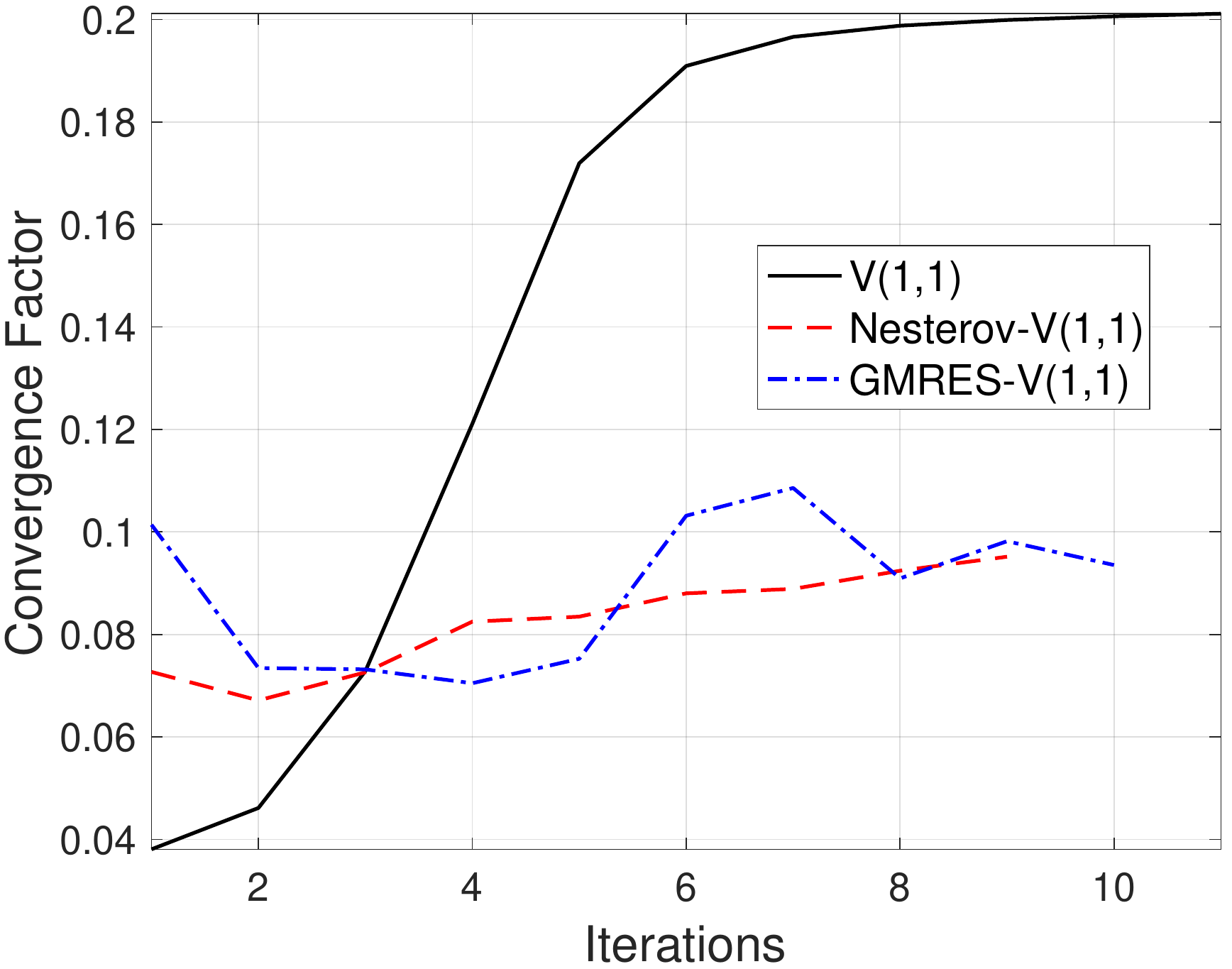}\label{fig:DiffusionExample_uniformConv}}
		\subfigure[]{\includegraphics[scale=0.3]{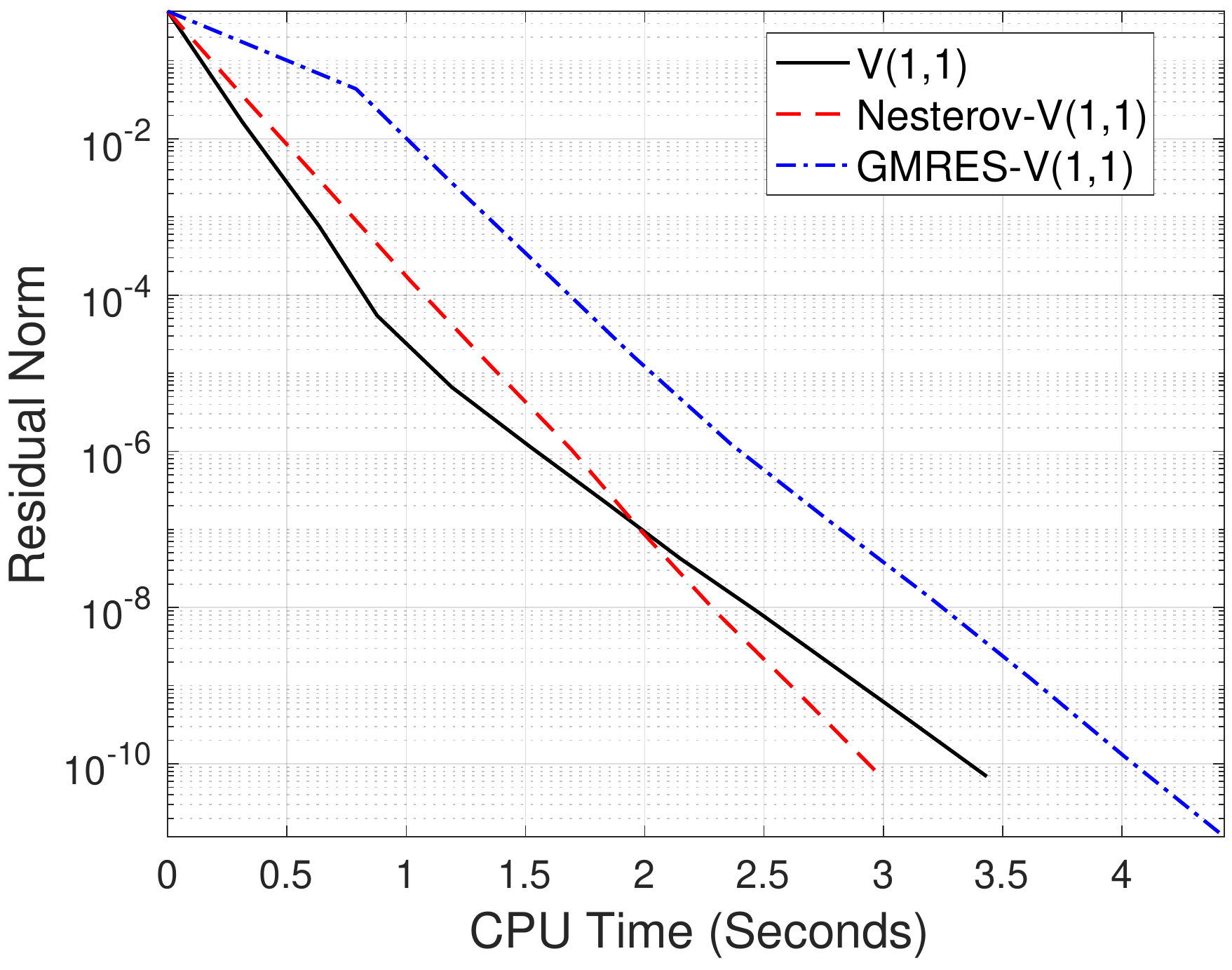}\label{fig:DiffusionExample_uniformCPU}}
		\caption{Comparison of accelerated Black Box Multigrid V$(1,1)$ cycles for the diffusion problem with log-normal (first row) and uniform (second row) distributions of the diffusion coefficient vector $\vsigma$.}
		\label{fig:DiffusionExample_lognormUniform}
	\end{figure}

	Finally, we test \eqref{eq:diffusionproblems}, \eqref{eq:DiscreteDiffusion} with a random diffusion coefficient vector $\vsigma$, sampled from a log-normal or uniform distribution, following \cite{greenfeld2019learning}. Due to the discontinuous coefficients, we use the classical Black Box Multigrid algorithm \cite{dendy1982black}, employing operator-dependent prolongation and Galerkin coarsening. For relaxation we employ Gauss-Seidel in natural (lexicographic) ordering. Following \cite{greenfeld2019learning}, we use bilinear finite element discretization. Here, we cannot use Fourier smoothing analysis, and $b_N$ is estimated by running V cycles with no acceleration, and we assume $b_1=0$. The results are shown in \Cref{fig:DiffusionExample_lognormUniform}. We find that Nesterov's scheme is competitive in terms of iteration count, and it is the fastest method in terms of CPU time in these experiments. 
	
	\section{Discussion and Conclusion}\label{sec:ConclusionandFuture}
	In this paper, we adapt Nesterov's scheme to accelerate stationary iterative methods for linear problems. Under the assumption that the eigenvalues of the iteration matrix are real, we derive a closed-form solution for the optimal scalar coefficient $c$ used in Nesterov's scheme. Numerical tests with accelerated multigrid cycles demonstrate the advantages of this approach. Moreover, we also study the robustness of Nesterov's scheme for cases where some of the eigenvalue of the iteration matrix $\mB$ are complex, identifying an explicit disk in the complex plane where the existence of complex eigenvalues does not degrade the rate of convergence. Our numerical results also demonstrate an advantage of Nesterov's acceleration scheme in such cases. 
	
	{\color{black} We note that the lower bound on the asymptotic convergence factor ACF shown in \cite{de2021asymptotic} is tight if the regime of the complex eigenvalues meets the disk defined in \Cref{them:Robustness}. Following the idea in \cite{de2021asymptotic}, the authors in \cite{wang2021asymptotic} quantified the improvement of the ACF  of Anderson acceleration  applied to Alternating Direction Method of Multipliers for nonlinear problems. From \cite{wang2021asymptotic}, we see that the Jacobians of ADMM at the fixed point always have complex eigenvalues and find that our results are still relevant to some nonlinear examples considered in \cite{wang2021asymptotic}. We also note that PCG and restricted-information Chebyshev are optimal when $\mB$ only has real eigenvalues and are faster than Nesterov's scheme. However, one may benefit from the multi-step acceleration described in \cite{de2021asymptotic} to reduce such a gap.} In future work we plan to study more general cases of complex eigenvalues not covered in this paper, {\color{black}to explore multi-step acceleration with optimal acceleration coefficients}, and to try to extend our approach to {\color{black} other methods for} nonlinear problems {\color{black} with analytic coefficients $c_k$}.

\section*{Acknowledgement}
We are grateful to Prof. Michael Zibulevsky for his advice. 

\appendix
\section{Proof of Lemma 1}\label{app:sec:proof:lemma:regionC}
	
We prove this lemma by showing that $r(c,b)\triangleq \max\{|\lambda_1(c,b)|,|\lambda_2(c,b)|\}>|b|$ for any nonzero $b$ if $|c|\geq 1$, so such a $c$ slows down convergence. 
	
	For $c > 1$, the first term in $\lambda_{1,2}$ already satisfies $\left|\frac{1}{2}(1+c)b\right| > |b|$ so either $\lambda_1$ or $\lambda_2$ (or both) must be larger than $b$. In the borderline case of $c=1$, the first term equals $b$, but in this case the second term cannot vanish except for $b=0$, which is excluded, so again either $\lambda_1$ or $\lambda_2$ must be strictly larger than $b$. 
	
	For $c \leq -1$ and $b>0$, $r(c,b)=-\lambda_2(c,b)$ and the square root term is real. To show that $-\lambda_2(c,b)>b$, we need to prove $\sqrt{(1+c)^2b^2-4cb}>(3+c)b$. Squaring both sides of this inequality, and simplifying, results in $(2+c)b<-c$, which is satisfied, because the left side is smaller than one in this regime, while the right side is at least one.
	
	For $c \leq -1$ and $b<0$, if $b > b_{cr}$ then the square root term is imaginary by \Cref{rem:bcr}, resulting by \Cref{rem:imaglambda} in $r(c,b)= \sqrt{cb} \geq \sqrt{-b} > |b|$. Otherwise, the square root is real and $r(c,b) = \lambda_1(c,b)$. To show that $\lambda_1(c,b)>|b|$, we need to prove that $\sqrt{(1+c)^2 b^2-4cb}>(-c-3)b$. This is indeed satisfied, due to the fact that the right side is negative because $-1 < b \leq b_{cr}$ implies $c < -3 - 2\sqrt{2}$.  \hfill \qedsymbol
\section{Proof of  Lemma 2}\label{app:sec:proof:lemma:extreme_rb}
	To prove that $r_c(b)$ has no local maximum for any $|c| < 1$, we first observe by \Cref{rem:bcr,rem:imaglambda} that in the range where $\lambda_{1,2}$ are complex,  $r_c(b) = \sqrt{cb}$ is strictly increasing (respectively, decreasing) for positive (respectively, negative) $b$, and therefore has no local maximum in this range. Outside this range, the derivative of $r_c(b)$ is discontinuous (only) at $b=0$ and $b = b_{cr}$. Elsewhere, it is given by 
	\begin{equation}
	r_c^\prime(b) = \frac{1}{2} \left( sgn(b)(1+c) +  \frac{(1+c)^2b-2c}{\sqrt{(1+c)^2b^2-4cb}} \right) .
	\label{eq:rcprime}
	\end{equation}
	
\begin{figure}
		\centering
		\subfigure[]{\includegraphics[scale = 0.45]{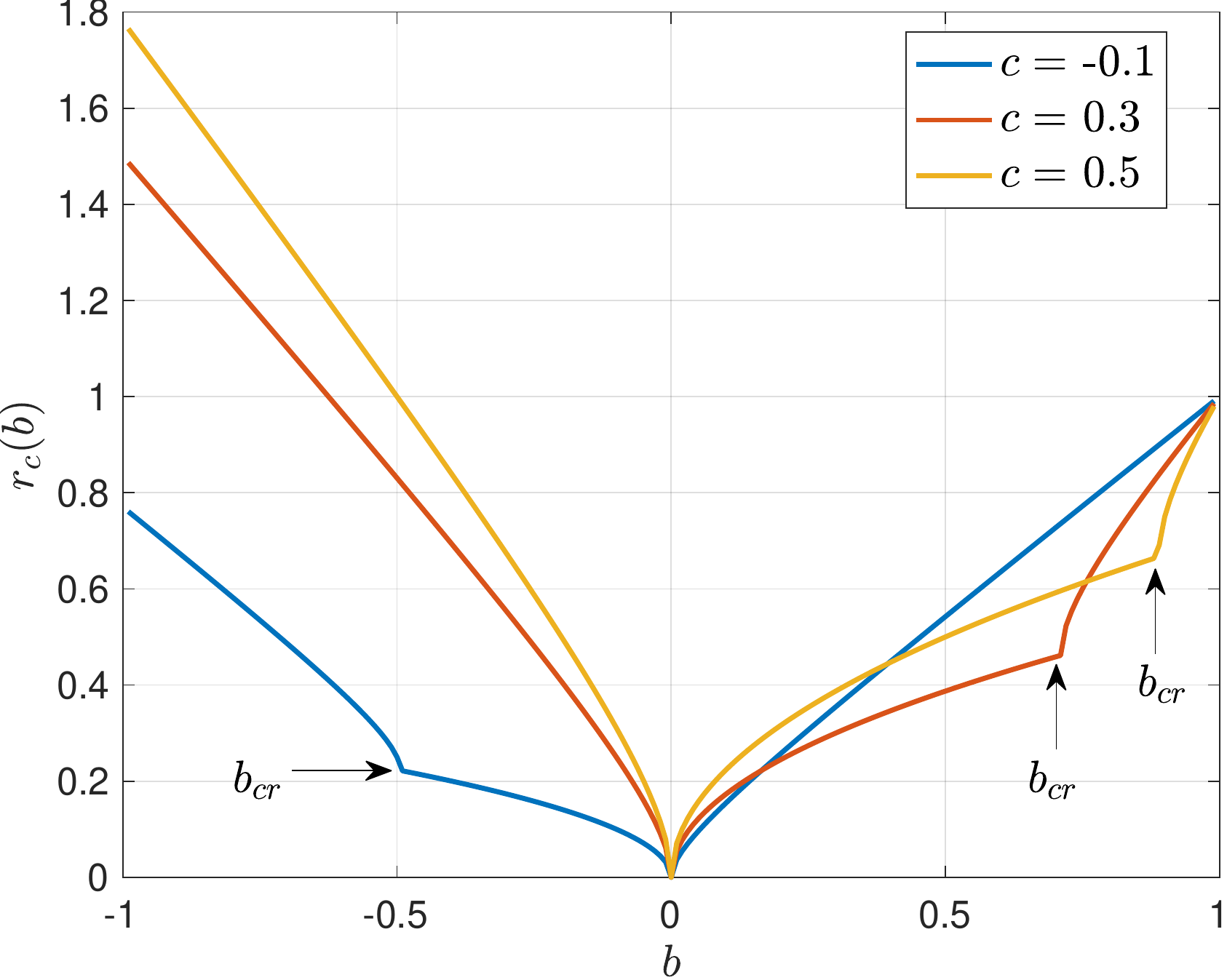}\label{fig:numericalshowrcb:r_c_b}}
		\subfigure[]{\includegraphics[scale = 0.45]{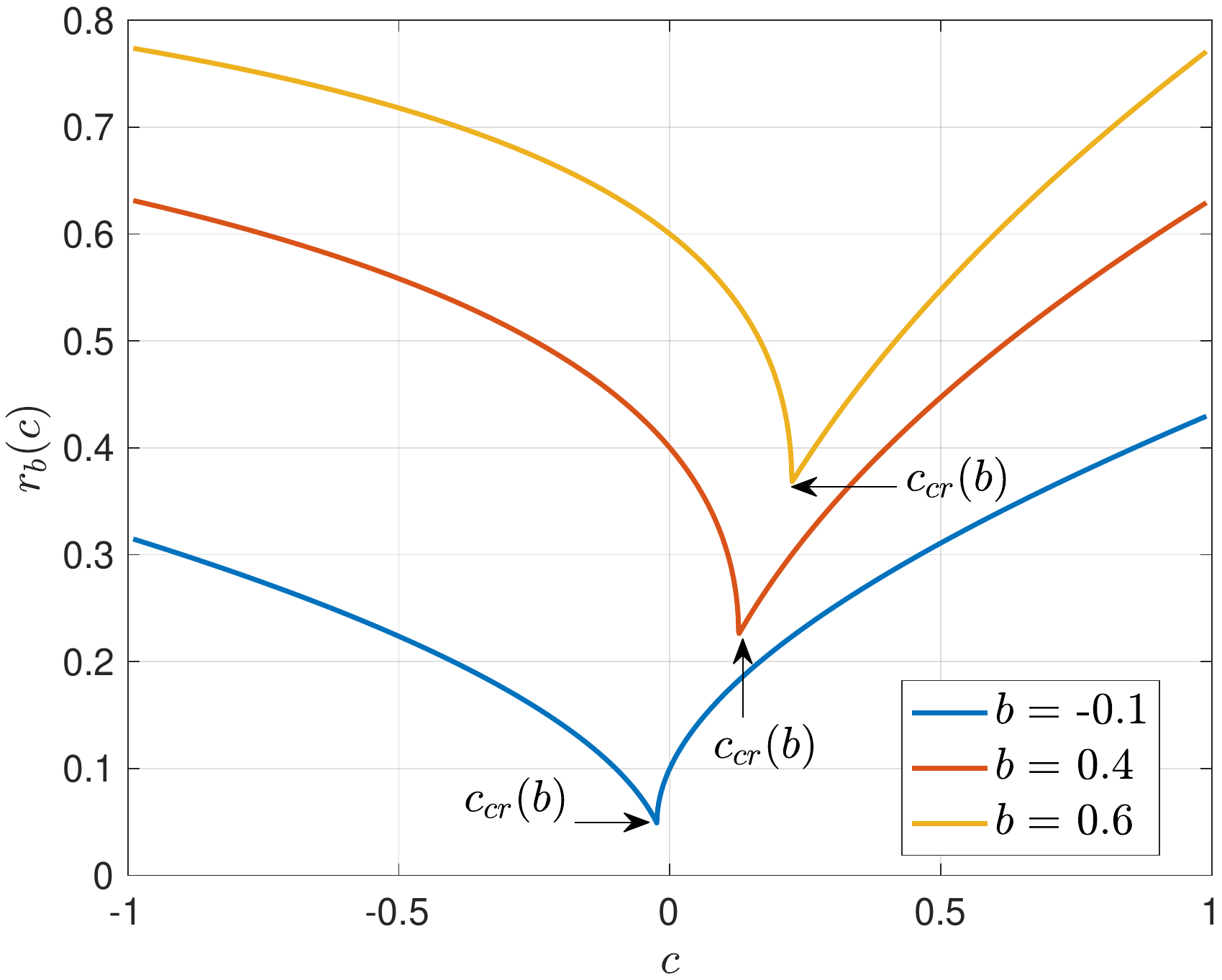}\label{fig:numericalshowrcb:r_b_c}}
		\caption{(a): The value of $r_c(b)$ as a function of $c$ and $b\in(-1,1)$. (b): The value of $r_b(c)$ as a function of $b$ and $c\in(-1,1)$.}\label{fig:numericalshowrcb}
	\end{figure}

	Note that for $b = b_{cr}$, the numerator in the second term is given by $2c$. It follows that $sgn(r_c^\prime(b_{cr}^-)) = sgn(r_c^\prime(b_{cr}^+))$ regardless of whether $c$ (hence also $b_{cr}$) is positive or negative, so $b_{cr}$ cannot be a local maximum. It remains to show that  $r^\prime_c(b)$ cannot vanish. We show this by comparing the squares of the first and second terms in the brackets on the right side of \eqref{eq:rcprime}:
	$$
	(1+c)^2 - \left( \frac{(1+c)^2b-2c}{\sqrt{(1+c)^2b^2-4cb}} \right)^2 = -\frac{4c^2}{(1+c)^2b^2-4cb} \, .
	$$
	This expression only vanishes for $c = 0$, for which \eqref{eq:rcprime} implies that $r^\prime_c(b) \neq 0$. Hence,  $r^\prime_c(b)$ cannot vanish for any $b \in (-1,1)$.  {\color{black}In \Cref{fig:numericalshowrcb}(a), we numerically show the value of $r_c(b)$ as a function of $c$ and $b\in(-1,1)$ that we clearly see the discontinuous at $b=0$ and $b=b_{cr}$.} 	

	\hfill \qedsymbol
\section{Proof of Lemma 3}\label{app:sec:proof:lemma:extreme_rc}
	Consider first the case $0 < b < 1$. In this regime, by \Cref{rem:ccr}, $r_b(c)=\sqrt{cb}$ for $c>c_{cr}(b)>0$, so $r^\prime_b(c) > 0$. On the other hand, for $c<c_{cr}(b)$  we have $r_b(c)=\frac{1}{2}\left((1+c)b+\sqrt{(1+c)^2b^2-4cb}\right)$, hence,
	\begin{eqnarray*}
	r^\prime_b(c) & = & \frac{1}{2} \left(b + \frac{(1+c)b^2 - 2b}{\sqrt{(1+c)^2b^2 - 4bc}} \right) \\ & = & \frac{b}{2} \left(1 - \sqrt{1+\frac{4(1-b)}{(1+c)^2b^2 - 4bc}} \right) < 0 \,.
	\end{eqnarray*}
	It follows that $r_b(c)$ is minimized in this regime for $c = c_{cr}(b)$. The derivation for $-1<b<0$ is analogous and we omit the details.
	{\color{black}In \Cref{fig:numericalshowrcb}(b), we numerically show the value of $r_b(c)$ as a function of $c$ and $b\in(-1,1)$ that $c_{cr}(b)$ is the minimizer.}
	
	\hfill \qedsymbol
\section{Proof of Thorem 1}\label{app:sec:proof:theorem:them:optimalC}
	
	In light of \Cref{lemma:extreme_rb,lemma:extreme_rc}, $c^*$ must satisfy one of the following three conditions: \circled{1} $c^* = c_{cr}(b_N)$, \circled{2} $r(c^*,b_1)=r(c^*,b_N)$, or  \circled{3} $c^* = c_{cr}(b_1)$, denoted $c_{top}$, $c_{mid}$ and $c_{bot}$, respectively. The choice is determined by the maximal $r$ obtained amongst the three. We assume for simplicity $|b_1|\leq b_N$, and remark that the complementary case, $|b_N|< -b_1$, is obtained analogously so the details are omitted.
	
	Note first that for $c<c_{cr}(b_N)$ we have $r_{c}(b_N) = \frac{1}{2} \left[(1+c)b_N +\sqrt{(1+c)^2b_N^2-4cb_N} \right]$, which is larger than $b_N$ if $c<0$, so it follows that $c^*\in[0,1)$. We first consider the case $b_1 \geq 0$. For $c<c_{cr}(b)$, we have for any $b$,
	$$
	r'_{c}(b)=\frac{1}{2}\left[(1 +c)+\frac{(1+c)^2b-2c}{\sqrt{(1+c)^2b^2-4cb}} \right]>0.
	$$
	On the other hand, for $c>c_{cr}(b)$, we have $r_{c}(b)=\sqrt{cb}$ so again $r'_{c}(b)>0$. The monotonicity of $r_c(b)$ implies $r_{c}(b_N)\geq r_{c}(b_1)$, with equality achieved only for $b_1=b_N$. It follows that $c^*=\argmin_c r_{b_N}(c)$ and case \circled{1} holds, i.e., $c^*=c_{top}\triangleq c_{cr}(b_N)$. The corresponding ACF is $r_{c_{top}}(b_N)=1-\sqrt{1-b_N}=\frac{2c_{top}}{1+c_{top}}$.
	
	It remains to consider the case $b_1<0$, whereby $r_{c}(b_1) = -\frac{1}{2} \left[ (1+c)b_1 - \sqrt{(1+c)^2b_1^2-4cb_1} \right]$. First, we identify the subrange in $b_1<0$ yielding $r_{c_{top}}(b_N)\geq r_{c_{top}}(b_1)$, leading us again to case \circled{1}. Denote for clarity $r_{top}=r_{c_{top}}(b_N).$ Then case \circled{1} holds if $b_1$ satisfies the following inequality,
	$$
	\begin{array}{rcl}
	r_{top}& \geq  &-\frac{1}{2} \left[(1+c_{top})b_1-\sqrt{(1+c_{top})^2b_1^2-4c_{top}b_1}\right], \\[8pt]
	2r_{top}+(1+c_{top})b_1&\geq &\sqrt{(1+c_{top})^2b_1^2-4c_{top}b_1}~~~(\text{square both sides and divide by 4}),\\
	r^2_{top}+r_{top}(1+c_{top})b_1&\geq &-c_{top}b_1,\\
	b_1&\geq &-\frac{r^2_{top}}{r_{top}(1+c_{top})+c_{top}}.
	\end{array}
	$$
	Using $r_{top}=\frac{2c_{top}}{1+c_{top}}$ and $b_N=\frac{4c_{top}}{(1+c_{top})^2}$, we obtain the condition $b_1\geq -\frac{1}{3}b_N$. 
	
\begin{figure}
		\centering
		\includegraphics[scale=0.45]{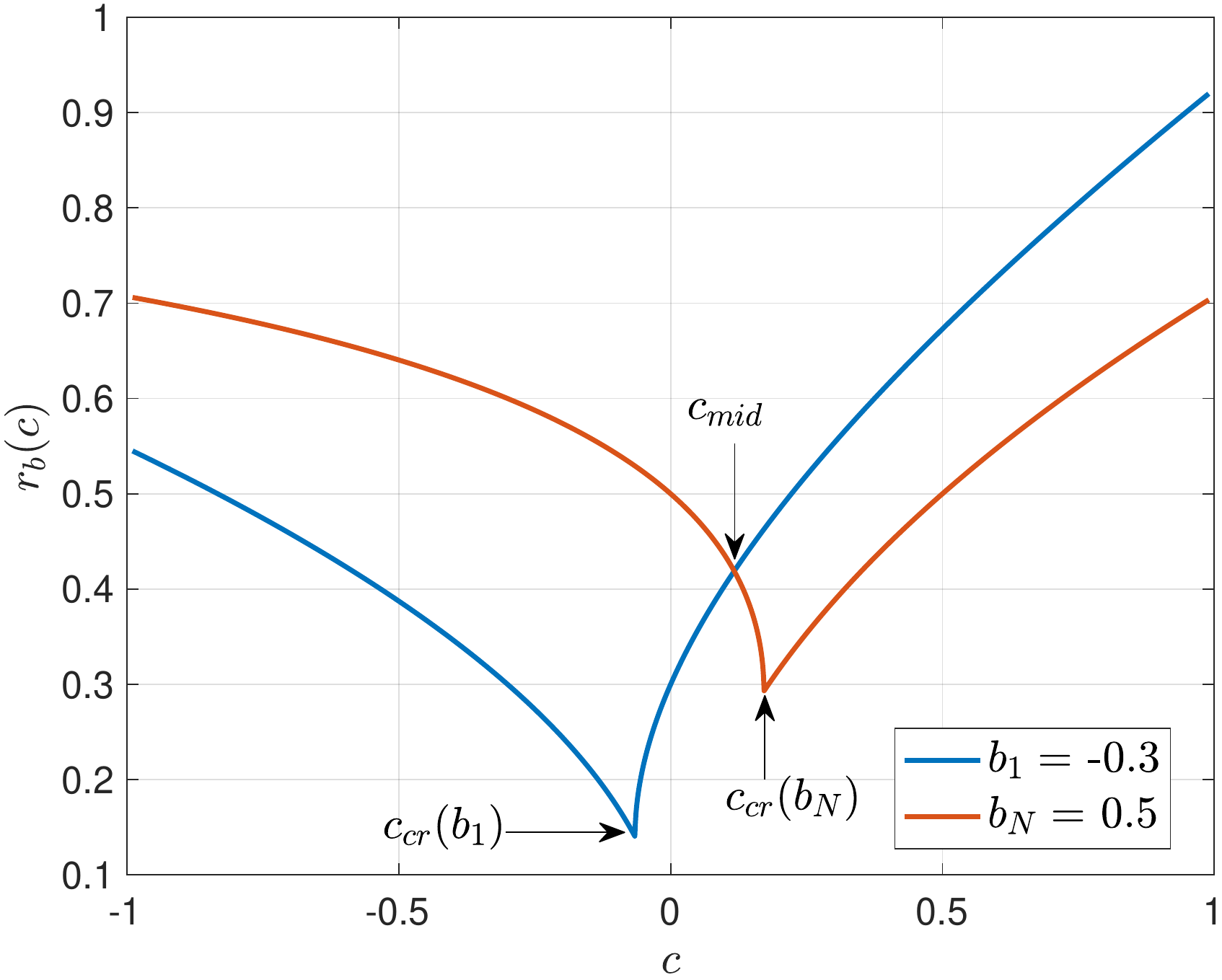}
		\caption{The curves of $r_b(c)$ as a function of $c$ with $b=b_1=-0.3$ and $b=b_N=0.5$.}\label{fig:Theorem_1_cross}
	\end{figure}

	Finally, we discuss the remaining range of $b_1<0$, i.e., $b_1\in[-b_N,-\frac{1}{3}b_N)$, where case \circled{1} does not hold. The fact that $c^*\in [0,1)$ rules out case \circled{3}, leaving case \circled{2}, that is,
	$$
	\begin{array}{rcl}
	-(1+c^*)b_1+\sqrt{(1+c^*)^2b_1^2-4c^*b_1}&=&(1+c^*)b_N+\sqrt{(1+c^*)^2b_N^2-4c^*b_N}.
	\end{array}\label{eq:r(c,b_1)=r(c,b_N)}
	$$
	Simplifying by repeatedly putting the square root terms on one side and squaring both sides, we get 
	$$
	\left[2b_1b_N\left(b_1+b_N\right)\right] (c^*)^2+\left[4b_1b_N\left(b_1+b_N\right)+\left(b_1-b_N\right)^2 \right]c^*+2b_1b_N\left(b_1+b_N\right) =0. \label{eq:b_1b_N}   
	$$
	Using $c^*\in[0,1)$, the valid solution of this quadratic equation is
	$$
	\begin{array}{rcl}
	c^*=c_{mid}&\triangleq&\frac{1-\sqrt{1-g(b_1,b_N)}}{1+\sqrt{1-g(b_1,b_N)}},
	\end{array}
	$$
	where $g(b_1,b_N) = - \frac{8b_1b_N(b_1+b_N)}{(b_1-b_N)^2}$.  {\color{black} In \Cref{fig:Theorem_1_cross}, we numerically show the curves of $r_b(c)$ as a function of $c$ and $b$ corresponding to \circled{2} that we clearly see $c_{mid}$ is the solution for $\min_c \max \{r_{b_1}(c),r_{b_N}(c)\}$.}

	As noted above, the derivation for $(b):~ |b_N|\leq -b_1$ is analogous and we omit the details.
\hfill \qedsymbol

\section{Proof of Theorem 2}\label{app:sec:proof:theorem:OptcComplexEigenvaluesB} 
	As in the proof of \Cref{them:optimalC},  we assume $|b_1|\leq b_N$ and remark that the complement, $~|b_N|< -b_1$, is proved analogously, omitting the details. 
	
	A complex eigenvalue $b^c$ does not influence the ACF, when $c^*$ is chosen according to \Cref{them:optimalC}, if $r(c^*,b^c)\leq r^*$, where $r^*=r_{c^*}(b_N)$, with $c^*=c_{top}$ or $c_{mid}$ as dictated by \Cref{them:optimalC}. This yields the explicit condition:
	\begin{equation}
	\frac{1}{2}\left|(1+c^*)b^c \pm \sqrt{(1+c^*)^2(b^c)^2-4c^*b^c}\right|\leq r^*.\label{eq:app:complexinequal}
	\end{equation}
	 Substituting $b^c = \bar b^c e^{j\theta}$ into \eqref{eq:app:complexinequal} and simplifying, we obtain
	\begin{equation}
	\frac{1+c^*}{2}\bar b^c\left|1 \pm \sqrt{1+\bar{b}e^{j\bar\theta}}\right|\leq r^*, \label{eq:app:complexinequalSimplify}
	\end{equation}
	where we have introduced the notation $\bar{b}=\frac{4c^*}{(1+c^*)^2\bar b^c}\geq0$ and $\bar\theta = \pi-\theta$. Since the term of the square root is complex and the sign of its real part is always positive, \eqref{eq:app:complexinequalSimplify} becomes
	\begin{equation}
	\frac{1+c^*}{2}\bar b^c\left|1 + \sqrt{1+\bar{b}e^{j\bar\theta}}\right|\leq r^*.\label{eq:app:complexinequalSimplifyFinal}
	\end{equation}
	Denote $\phi(\bar\theta) = \left|1 + \sqrt{1+\bar{b}e^{j\bar\theta}}\right|$. To bound the left-hand side from above, we derive a monotonicity property for $\phi(\bar\theta)$ in the period $\bar\theta \in (-\pi,\pi]$. With Euler's formula, we can rewrite $\phi(\bar\theta)$ as $\sqrt{1+2\chi\cos\frac{\vartheta}{2}+\chi^2}$, where $\chi=\left(1+2\bar{b}\cos\bar\theta+\bar{b}^2\right)^{\frac{1}{4}}$ and  $\vartheta \in(-\frac{\pi}{2},\frac{\pi}{2}): =\arctan \frac{\bar{b}\sin{\bar\theta}}{1+\bar{b}\cos{\bar\theta}}$, and then it is evident that $\phi(-\bar\theta) = \phi(\bar\theta)$. Since $\phi(\bar\theta)$ is an even function, we consider its monotonicity in the half period $\bar\theta \in (0,\pi)$.  Using $\sqrt{x+jy}=\sqrt{\frac{x+\sqrt{x^2+y^2}}{2}}+sgn(y)j\sqrt{\frac{-x+\sqrt{x^2+y^2}}{2}}$ and $\phi(\bar\theta)=\left|1+\sqrt{1+\bar{b}\cos\bar\theta+j\bar{b}\sin\bar\theta}\right|$, we get $\phi(\bar\theta)=\left|1+\sqrt{\frac{1+\bar{b}\cos\bar\theta+\chi^2}{2}}+j\sqrt{\frac{-(1+\bar{b}\cos\bar\theta)+\chi^2}{2}}\right|$. For convenience, we consider the monotonicity of $\phi^2(\bar\theta)=1+\sqrt{2+2\bar b\cos\bar\theta+2\chi^2}+\chi^2$ instead of $\phi(\bar\theta)$. Since 
	$$
	\frac{\partial \phi^2(\bar\theta)}{\partial \bar\theta}= 2\chi \frac{\partial \chi}{\partial\bar\theta}+\frac{-\bar{b}\sin\bar\theta+2\chi\frac{\partial \chi}{\partial \bar\theta}}{\sqrt{2(1+\bar{b}\cos\bar\theta)+2\chi^2}}< 0,~~~\bar\theta\in (0,\pi)
	$$
	(using $\frac{\partial \chi}{\partial\bar\theta}<0$, $\bar{b}\sin\bar\theta\geq0$, and $\chi\geq 0$), we find  that $\phi(\bar\theta)$ is monotonically decreasing in $\bar\theta \in (0,\pi)$. Together with the fact that $\phi(\theta)$ is an even function, we have $r(c^*,\bar b^c)\leq r(c^*,b^c)\leq r(c^*,-\bar b^c)$.
	
	From \Cref{them:optimalC}, we know that $r^{*}=r_{c_{top}}(b_N)\geq r_{c_{top}}(b_i)$ for any real $b_i\geq -\frac{1}{3}b_N$ when $T_{top}$ is the relevant domain. Thus, for any complex eigenvalue with $\bar b^c\leq \frac{1}{3}b_N$, we obtain $r(c_{top},b^c)\leq r^*$. Since $\phi(\theta)$ is monotonically increasing (respectively, decreasing) in $\theta\in(0,\pi)$ (respectively, $\theta\in(-\pi,0)$), the upper bound on $r(c_{top},b^c)$ becomes tight when the argument is close to $\pi$ (respectively, $-\pi$). This means that if the modulus is larger than $\frac{1}{3}b_N$ then the argument should be close to $0$ to enforce $r(c_{top},b^c)\leq r^*$. Evidently, the conclusion drawn here is valid also for $c^*=c_{mid}$, as it is implied by the monotonicity properties of $\phi(\theta)$. 
	
	As noted above, the proof for the complementary domain $|b_N|< -b_1$ is analogous and so the details are omitted.  \hfill \qedsymbol

\bibliographystyle{unsrt}  
\bibliography{refs}



\end{document}